\documentclass[11pt]{amsart}
\usepackage{amsmath,amssymb,amsfonts,url}
\numberwithin{equation}{section}

\newtheorem{thm}{Theorem}[section]
\newtheorem{lem}{Lemma}[section]
\newtheorem{rem}{Remark}[section]
\newtheorem{prop}{Proposition}[section]

\newtheorem{mainthm}{Theorem}

\newcommand{\hdot}{^\text{\r{}}\hspace{-.33cm}H}

\begin{document}
\title[Liouville system]{A topological degree counting for some Liouville systems of Mean Field Equations} 

\author{Chang-shou Lin}
\address{Department of Mathematics\\
        Taida Institute of Mathematical Sciences\\
        National Taiwan University\\
         Taipei 106, Taiwan } \email{cslin@math.ntu.edu.tw}

\author{Lei Zhang}
\address{Department of Mathematics\\
        University of Florida\\
        358 Little Hall P.O.Box 118105\\
        Gainesville FL 32611-8105}
\email{leizhang@math.ufl.edu}
\thanks{Zhang is supported in part by NSF Grant 0900864 (1027628)}

\date{\today}

\begin{abstract}
Let $A=(a_{ij})_{n\times n}$ be an invertible matrix and $A^{-1}=(a^{ij})_{n\times n}$ be the inverse of $A$.
In this paper, we consider the generalized Liouville system:
\begin{equation}\label{abeq1}
\Delta_g u_i+\sum_{j=1}^n a_{ij}\rho_j\left(\frac{h_j e^{u_j}}{\int h_j e^{u_j}}-1\right)=0\quad\text{in \,}M,
\end{equation}
where $0< h_j\in C^1(M)$ and $\rho_j\in \mathbb R^+$, and prove that, under the assumptions of $(H_1)$ and $(H_2)$\,(see Introduction),
the Leray-Schauder degree of \eqref{abeq1} is equal to
\begin{equation*}
\frac{(-\chi(M)+1)\cdots (-\chi(M)+N)}{N!}
\end{equation*}
if $\rho=(\rho_1,\cdots, \rho_n)$ satisfies
\begin{equation*}
8\pi N\sum_{i=1}^n\rho_i<\sum_{1\leq i,j\leq n}a_{ij}\rho_i\rho_j<8\pi(N+1)\sum_{i=1}^n\rho_i.
\end{equation*}
Equation \eqref{abeq1} is a natural generalization of the classic Liouville equation and is the Euler-Lagrangian equation of Nonlinear function $\varPhi_\rho$:
\begin{equation*}
\varPhi_\rho(u)=\frac{1}{2}\int_M\sum_{1\leq i,j\leq n}a^{ij}\nabla_g u_i\cdot \nabla_g u_j+\sum_{i=1}^n\int_M\rho_iu_i
-\sum_{i=1}^n\rho_i\log \int_M h_i e^{u_i}.
\end{equation*}
The Liouville system \eqref{abeq1} has arisen in many different research areas in mathematics and physics.
Our counting formulas are the first result in degree theory for Liouville systems.
\end{abstract}

\maketitle

\section{Introduction}

Let $(M,g)$ be a compact Riemann surface with volume $1$, $h_1,...,h_n$ be positive $C^1$ functions on $M$, $\rho_1,..,\rho_n$ be nonnegative
constants. In this article we consider the following Liouville system defined on $(M,g)$:
\begin{equation}\label{mainsys}
\Delta_g u_i+\sum_{j=1}^n\rho_j a_{ij} (\frac{h_je^{u_j}}{\int_Mh_je^{u_j}dV_g}-1)=0,\quad i\in I:=\{1,..,n\}
\end{equation}
where $dV_g$ is the volume form, $\Delta_g$ is the Laplace-Beltrami operator, in local coordinates it is of the form:
$$\Delta_g=\sum_{i,j=1}^2\frac{1}{\sqrt{\mbox{det}(g)}}\frac{\partial}{\partial x_i}(\sqrt{\mbox{det}(g)}g^{ij}\frac{\partial}{\partial x_j}),\quad
(g^{ij})_{2\times 2}=(g_{ij})_{2\times 2}^{-1}.$$

\noindent
When $n=1$, equation \eqref{mainsys} is the mean field equation of the Liouville type:
\begin{equation}\label{equfromsys}
\Delta_g u+\rho\left(\frac{h e^u}{\int_M h e^udV_g}-1\right)=0\quad\text{in \,}M
\end{equation}
when $a_{11}=1$. Therefore, the Liouville system \eqref{mainsys} is a natural extension of the classical Liouville equation,
which has profound connection with geometry and physics, and has been extensively studied for the past three decades.

If $u$ is a solution of \eqref{mainsys}, then after adding a constant, $u+c$ is also a solution of \eqref{mainsys}.
Hence, we can always assume $u\in \,\hdot^1(M)$, where
\begin{equation*}
\hdot^1(M)=\left\{
u\in L^2(M)\,\Big|\, |\nabla_g u|\in L^2(M)\text{ and }\int_M  u\,dV_g=0
\right\}.
\end{equation*}
For any $\rho=(\rho_1,\cdots, \rho_n)$, $\rho_i>0$, let $\varPhi_\rho$ be a nonlinear functional defined in \,\,$\hdot^{1,n}=\,\hdot^1(M)\times\cdots\times\,\hdot^1(M)$ by
$$
\varPhi_\rho(u)=\frac{1}{2}\sum_{i,j\in I} a^{ij}\int_M \nabla_g u_i\cdot \nabla_g u_j dV_g-\sum_{j\in I}\rho_j\log\int_M h_j e^{u_j}dV_g
$$
where $(a^{ij})_{n\times n}$ is the inverse of $A=(a_{ij})_{n\times n}$.
It is easy to see that equation \eqref{mainsys} is the Euler-Lagrangian equation of $\varPhi_\rho$.

For a bounded smooth domain $\Omega$ in $\mathbb R^2$, we are also interested in the following system of equations:
\begin{equation}\label{systemtrans}
\begin{cases}
&\Delta u_i+\sum_{j=1}^n a_{ij}\rho_j\frac{h_j e^{u_j}}{\int_\Omega h_je^{u_j}dx}=0\quad\text{in \,}\Omega,\\
&u_i|_{\partial\Omega}=0, \quad i\in I
\end{cases}
\end{equation}
where $h_j$ are positive $C^1$ function on $\overline \Omega$.

The Liouville equation \eqref{equfromsys} or systems \eqref{mainsys} and \eqref{systemtrans} have appeared in many different disciplines in mathematics.
In conformal geometry, when $\rho=8\pi$ and $M$ is the sphere $\mathbb S^2$, equation \eqref{equfromsys} is equivalent to the famous Nirenberg problem.
For a bounded domain in $\mathbb R^2$ and $n=1$, \eqref{systemtrans} can be derived from the mean field limit of Euler flows or spherical Onsager vortex theory,
as studied by Caglioti, Lions, Marchioro and Pulvirenti\,\cite{lion1, lion2}, and Kiessling\,\cite{kiess}, Chanillo and Kiessling\,\cite{chanillo1}.
In classical gauge field theory, equation \eqref{mainsys} is closely related to the Chern-Simons-Higgs equation for the abelian case,
see \cite{caffarelli,hong,jackiw,yang}.
Various Liouville systems are also used to describe models in the theory of Chemotaxis\,\cite{childress, keller}, in the physics of charged particle beams\,\cite{bennet,debye,kiess2,kiessling2},
in the non-abelian Chern-Simons-Higgs theory\,\cite{dunne,jostlinwang,yang} and other gauge field models \cite{dziar, kimleelee}.
For recent developments of these subjects or related Liouville systems in more general
settings, we refer the readers to \cite{aly, barto1, barto3, biler, ChenLin1, chenlin2, chenlinnew, clwang, CSW, CSW1, licmp, linduke, linarch, linzhang1, nolasco2, rubinstein, spruck, wolansky1, wolansky2, zhangcmp, zhangccm} and the references therein.

For a bounded smooth domain $\Omega$ in $\mathbb R^2$,
Chipot-Shafrir-Wolansky\,\cite{CSW} considered equation \eqref{systemtrans}, where the constant matrix $A=(a_{ij})_{n\times n}$ satisfies the following condition:

\medskip

\noindent
{\bf (H1)}: $A$ is symmetric, nonnegative, irreducible and invertible.

\medskip

\noindent
Here $A$ is called nonnegative if $a_{ij}\geq 0$ for all $i,j\in I=\{1,2,\cdots,n\}$, and is called irreducible if there is no subset $J$ of $I$ such that
\begin{equation*}
a_{ij}=0\quad\text{for all }i\in J\text{ and }j\in I\setminus J.
\end{equation*}

\noindent
In another word, equation \eqref{systemtrans} can not be written as two de-coupled subsystems.\,In \cite{CSW}, the authors introduced
nonlinear functions $\Lambda_J(\rho)$ of $\rho=(\rho_1,\cdots, \rho_n)$ defined by
\begin{equation*}
\Lambda_J(\rho)=8\pi\sum_{i\in J}\rho_i-\sum_{i,j\in J}a_{ij}\rho_i\rho_j
\end{equation*}
for any non-empty subset $J$ of $I=\{1,2,\cdots, n\}$. Let
\begin{equation*}
\Gamma=\{
(\rho_1,\cdots,\rho_n)\,|\,\rho_i>0, \, \Lambda_J(\rho)>0\text{ for all }\emptyset\subsetneq J\subsetneq I
\}.
\end{equation*}
\noindent
Among other things, Chipot-Shafrir-Wolansky\,\cite{CSW} proved the following theorem.
\begin{mainthm}\label{mainathm}
Suppose $A$ satisfies {\bf (H1)}, $h_1,..,h_n$ are positive $C^1$ functions on $\overline\Omega$, and $\rho=(\rho_1,\cdots,\rho_n)$ satisfies
\begin{equation}\label{sufcondlamj}
\rho\in\Gamma.
\end{equation}
Then equation \eqref{systemtrans} possesses a solution.
\end{mainthm}

We note that in \cite{CSW,CSW1}, the authors also proved that the sufficient condition \eqref{sufcondlamj} in Theorem \ref{mainathm}
is also a necessary condition for the existence of equation \eqref{systemtrans} when $\Omega$ is a ball.
When $n=1$, equation \eqref{systemtrans} with the parameter $\rho\in\Gamma$ is equivalent to the Liouville equation with $\rho<8\pi$.
From various known results of Liouville equations, we expect that solutions of equation \eqref{systemtrans} with $\rho\not\in\Gamma$, should have Morse index bigger than $1$.
Therefore, the classical Leray-Schauder degree theory is a suitable tool to be applied for studying equation \eqref{mainsys} or \eqref{systemtrans} when $\rho\not\in\Gamma$.

To apply the degree theory, we should first prove the a priori bound for non-critical parameter $\rho$,
or equivalently, study the asymptotic behavior of bubbling solutions.
In \cite{linzhang1}, we have proved that near each blow up point, the behavior of these bubbling solutions can be controlled well
by the standard bubble, under the assumption that all the components $u_i$ of solutions have to blow up simultaneously i.e.,
a suitable scaling of $u_i$ should converge to some entire solutions of Liouville system:
\begin{equation}\label{newentiresol}
\begin{cases}
&\Delta U_i+\sum_{j\in I}a_{ij} e^{U_j}=0\quad\text{in \,}\mathbb R^2,\\
&\int_{\mathbb R^2} e^{U_j}dx<+\infty,\quad i\in I.
\end{cases}
\end{equation}
This blow-up is called type 1. However, there might be situations that at some blow-up points,
only part of the components of $u_i$ ($i\in I$) blows up, but the remaining part does not blow up.
This blow-up is called type 2. If both type 1 and type 2 occur simultaneously for a sequence of bubbling solution, then the set of critical parameters could be very complicated.
For example $n=2$, the critical parameter might be
\begin{align*}
&\big\{
(\rho_1,\rho_2)\,|\,\rho_1=\tilde \rho_1+8\pi m, \rho_2=\tilde\rho_2+8\pi l\text{ and }(\tilde\rho_1,\tilde\rho_2)\text{ satisfies }\\
&8\pi k(\tilde\rho_1+\tilde\rho_2)=\sum_{i,j}a_{ij}\tilde\rho_1\tilde\rho_2, \,m.l,k\in\mathbb N
\big\},
\end{align*}
where $a_{11}=a_{22}=1$ is assumed, $\mathbb N$ is the set of all natural numbers. Thus, the topological degree would be very difficult to compute.
In this paper, we will prove that this complexity can be avoided if we assume the coefficient matrix $A$ to satisfy {\bf (H1)} and the following {\bf (H2)} condition:

\medskip

\noindent
{\bf (H2)}: $a^{ii}\leq 0$ for $i\in I$, $a^{ij}\geq 0$ for $i\ne j$, $i,j\in I$,\\
\quad\quad\quad and $\sum_{j\in I}a^{ij}\geq 0$ for $i\in I$.

\medskip

\noindent
Throughout the paper, we assume that $A$ satisfies both {\bf (H1)} and {\bf (H2)}. For $n=2$, $A=\begin{pmatrix}a_{11}&a_{12}\\a_{12}&a_{22}\end{pmatrix}$ satisfies {\bf (H1)} and {\bf (H2)}
if and only if $a_{ij}\ge 0$, $\max(a_{11},a_{22})\leq a_{12}^2$ and $\det A\ne 0$. For $n=3$, assume $A=\begin{pmatrix}0& a_1 &a_2\\a_1& 0& a_3\\ a_2& a_3&0\end{pmatrix}$.
Then $A$ satisfies {\bf (H1)} and {\bf (H2)} if and only if $a_i>0$ and $a_i+a_j\geq a_k$ for $i\ne j\ne k$.

To state our results, we begin with equation \eqref{mainsys}.
\begin{thm}\label{main1}
Let $A=(a_{ij})_{n\times n}$ satisfy {\bf (H1)} and {\bf (H2)}, and $N$ be a nonnegative integer and
\begin{align*}
\mathcal{O}_N&=\big\{
(\rho_1,\cdots,\rho_n)\,\big|\, \rho_i\geq 0, i\in I\text{ and }\\
&8\pi N\sum_{i\in I}\rho_i<\sum_{i,j\in I}a_{ij}\rho_i\rho_j<8\pi (N+1)\sum_{i\in I}\rho_i
\big\}.
\end{align*}
Suppose $h_i\in C^1(M)$ is positive and $K$ is a compact subset of $\mathcal{O}_N$. Then there exists a constant $C$ such that for any solution $u=(u_1,\cdots, u_n)$ of \eqref{mainsys} with $\rho\in K$, we have
\begin{equation*}
|u_i(x)|\leq C\quad\text{for \,}i\in I\text{ and }x\in M.
\end{equation*}
\end{thm}

Note that the set $\mathcal{O}_N$ is bounded if $a_{ii}>0$ for all $i$, and $\mathcal{O}_N$ is unbounded if $a_{ii}=0$ for some $i$. By Theorem \ref{main1},
the critical parameter set for \eqref{mainsys} is the set $\Gamma_N$, where
\begin{equation*}
\Gamma_N=
\left\{
\rho\,\big|\, 8\pi N \sum_{i\in I}\rho_i=\sum_{i,j\in I}a_{ij}\rho_i\rho_j
\right\}.
\end{equation*}

After Theorem \ref{main1}, for $\rho\not\in \Gamma_N$ for any positive integer $N$, we can define the nonlinear map $T_\rho=(T^1,\cdots, T^n)$ from \,\,$\hdot^{1,n}=\,\hdot^1(M)\times\cdots\times\,\hdot^1(M)$
to \,$\hdot^{1,n}$ by
\begin{equation*}
T^i=-\Delta^{-1}_g\left(\sum_{j\in I} a_{ij}\rho_j\left(\frac{h_j e^{u_j}}{\int_M h_j e^{u_j}}-1\right)\right),\quad i\in I.
\end{equation*}
Obviously, $T_\rho$ is compact from \,$\hdot^{1,n}$ to itself. Then we can define the Leray-Schauder degree of equation \eqref{mainsys} by
\begin{equation*}
d_\rho={\rm deg}\, (I-T_\rho\,;\, B_R,0),
\end{equation*}
where $R$ is sufficiently large and $B_R=\{u\,|\,u\in \,\hdot^{1,N}\text{ and }\|u\|<R\}$. By the homotopic invariance and Theorem \ref{main1}, $d_\rho$ satisfies the following properties:

\medskip

\noindent
(i) $d_\rho$ is a constant for $\rho\in \mathcal{O}_N$,

\medskip

\noindent
(ii) $d_\rho$ is independent of $h=(h_1,h_2,\cdots, h_n)$.

\medskip

\noindent
The following result is the formula for computing $d_\rho$.
\begin{thm}\label{main2}
$d_{\rho}$ be the Leray-Schauder degree for (\ref{mainsys}) for $\rho\in \mathcal{O}_{N-1}$, $N\in \mathbb N$. Then
\begin{equation}\label{dcount}
d_{\rho}=\left\{\begin{array}{ll}
1\quad \text{if }\rho \in \mathcal{O}_0\\
\frac{1}{N!}\bigg ((-\chi_M+1)...(-\chi_M+N) \bigg )\quad \text{if }\rho \in \mathcal{O}_N.
\end{array}
\right.
\end{equation}
where $\chi_M$ is the Euler characteristic of $M$.
\end{thm}
Since $\chi_M=2-2g_e$ where $g_e$ is the genus of $M$, the following existence theorem is
implied by Theorem \ref{main2}:

\begin{thm}\label{exist1}(Main Theorem)
Let $M$ be a compact Riemann surface with genus greater than $0$ and $h_1,..,h_n$ be positive $C^1$ functions on $M$. Then (\ref{mainsys}) always has a solution for $\rho\not\in\Gamma_N$ for any $N\in \mathbb N$.
\end{thm}

Similarly, for equation \eqref{systemtrans}, we have the following result:

\begin{thm}\label{dirthm}
Let $(h_1,..,h_n)$ be positive $C^1$ functions on $\bar \Omega$. Then the Leray-Schauder degree $d_{\rho}$ for \eqref{systemtrans} is
$$d_{\rho}
=\left\{\begin{array}{ll}
1,\quad \rho \in \mathcal{O}_0\\
\frac{1}{N!}\bigg ((-\chi+1)...(-\chi+N) \bigg ),\, \rho \in \mathcal{O}_N,\quad N\in \mathbb N
\end{array}
\right.
$$
where $\chi=1-g_e$, $g_e$ is the number of holes inside $\Omega$. In particular if $\Omega$ is not simply connected, \eqref{systemtrans} always has a solution for $\rho$.
\end{thm}

The organization of this paper is as follows.
In section 2 we mainly address entire solutions of \eqref{newentiresol} and some important properties implied by the assumption {\bf (H2)}.
Then in section 3 we give a detailed description on the asymptotic behavior of blowup solutions near a blowup point.
In \cite{linzhang1} the authors have proved that if the system all converges to an entire system of $n$ equations around each blow-up,
then all the blow-up solutions converge to the same system after scaling.
In this section we consider the case of the type 2 blow-up, and give a sharper estimate for the bubbling solution near the blow-up point.
In section 4 by using the sharper estimates in section 3, we will prove that type 1 and type 2 blow-up can not occur simultaneously in any sequence of bubbling solutions,
which leads to the proof of Theorem \ref{main1} and Theorem \ref{main2} in section 5. Finally in section 6 we prove Theorem \ref{dirthm}.

\medskip

\noindent{\bf Acknowledgement}
Part of the paper was finished when the second author was visiting Taida Institute for Mathematical Sciences (TIMS) in May 2010. He is very grateful to TIMS for their warm hospitality. He also would like to thank the National Science Foundation for the partial support (NSF-DMS-10276280).

\section{On entire solutions}

In this section, we discuss the entire solutions of the Liouville system
\begin{equation}\label{oureq}
\left\{\begin{array}{ll}
\Delta u_i+\sum_{j\in I}a_{ij}e^{u_j}=0,\quad \mathbb R^2,\quad i\in I\\
\\
\int_{\mathbb R^2}e^{u_i}<\infty,\quad i\in I.
\end{array}
\right.
\end{equation}
System (\ref{oureq}) is closely related to the following system of equations
\begin{equation}\label{chipotshi}
\left\{\begin{array}{ll}
\Delta v_i+\mu_i e^{\sum_{j\in I}a_{ij}e^{v_j}}=0\quad \mathbb R^2,\quad i\in I. \\
\\
\int_{\mathbb R^2}\mu_i e^{\sum_{j\in I}a_{ij}e^{v_j}}<\infty,\quad i\in I,
\end{array}
\right.
\end{equation}
which was studied initially by Chanillo and Kiessling \cite{chanillo2}, and
by Chipot-Shafrir-Wolansky\,\cite{CSW,CSW1}. Obviously, these two systems are equivalent if the coefficient matrix $A=(a_{ij})$ is invertible.
In this section, the coefficient matrix $A$ is assumed to be symmetric and nonnegative, but not necessarily invertible. After a permutation of rows, $A$ can be written as
\begin{equation}\label{matrixa}
\left(\begin{array}{cccc}
A_1&  &  & \\
  & A_2 &  & \\
  &           & \ddots &  \\
  &   &   & A_k
  \end{array}
  \right )
\end{equation}
where each $A_l:=(a_{ij})_{i,j\in I_l}$\,$(l=1,\cdots,k)$ is irreducible and $I=\bigcup_{l=1}^k I_l$.
For a positive vector $\sigma=(\sigma_1,\cdots,\sigma_n)$\,(which means each $\sigma_i$ is positive), we define
\begin{equation*}
\Lambda_J(\sigma)=4\sum_{i\in J}\sigma_i-\sum_{i,j\in J}a_{ij}\sigma_i\sigma_j\quad\text{for any \,}\emptyset\subsetneq J \subset I.
\end{equation*}
Then Chipot-Shafrir-Wolansky prove the following theorem in \cite{CSW}:

\medskip

\noindent
\emph{{\rm {\bf Theorem B.}} (Theorem 1.4 of \cite{CSW}): Let $A$ be a nonnegative, symmetric matrix that satisfies (\ref{matrixa}). Let
$\sigma=(\sigma_1,\cdots,\sigma_n)$ be a positive vector such that
\begin{equation}\label{iffcon}
\Lambda_{I_l}(\sigma)=0,\quad \Lambda_J(\sigma)>0,\quad \emptyset\subsetneq J\subsetneq I_l,\quad l=1,\cdots,k.
\end{equation}
Then there exist $v=(v_1,\cdots,v_n)$ and $\mu=(\mu_1,\cdots,\mu_n)$ ($\mu_i>0$) such that (\ref{chipotshi}) is satisfied and
\begin{equation}\label{519e1}
\sigma_i=\frac 1{2\pi}\int_{\mathbb R^2}\mu_ie^{\sum_{j\in I}a_{ij}v_j}, \quad i\in I.
\end{equation}
Conversely, for an entire solution $v$ to (\ref{chipotshi}), $\sigma$ defined by (\ref{519e1}) satisfies (\ref{iffcon}).
}
\medskip

We note that if one submatrix $A_l$ is zero, then $I_l$ consists of only one element because $A_l$ is irreducible. In this case, no positive
$\sigma$ satisfies $\Lambda_{I_l}(\sigma)=4\sigma=0$. Therefore (\ref{chipotshi}) has no solution for any positive $\mu=(\mu_1,\cdots,\mu_n)$.

\medskip

Let $\sigma_i$ ($i\in I$) be positive such that $\sigma=(\sigma_1,\cdots,\sigma_n)$ satisfies \eqref{iffcon}.
Then by Theorem B, there is a solution $v=(v_1,\cdots, v_n)$ of \eqref{chipotshi} such that \eqref{519e1} holds.
In \cite{CSW}, Chipot-Shafrir-Wolansky proved the asymptotic behavior of $v$ at $\infty$\,:
\begin{equation}
v_i(x)=-\sigma_i\log |x|+O(1),\quad i\in I,
\end{equation}
and
\begin{equation}\label{new27}
\sum_{j\in I}a_{ij}\sigma_j>2,
\end{equation}
due to $\exp\left(\sum_{j\in I}a_{ij}v_j\right)\in L^1(\mathbb R^n)$. For any solution $v$ of \eqref{chipotshi}, set
\begin{equation}
u_i=\sum_{j\in I}a_{ij}v_j+\log \mu_i.
\end{equation}
Then $u_i$ is a solution of the system \eqref{oureq} with
\begin{equation}\label{new210}
u_i=-m_i \log|x|+O(1)\quad\text{for large \,}|x|,
\end{equation}
and
\begin{equation}
m_i>2,
\end{equation}
where
\begin{equation}\label{new28}
m_i=\sum_{j\in I}a_{ij}\sigma_j.
\end{equation}

\medskip

Conversely, let $u=(u_1,\cdots, u_n)$ be a solution of \eqref{oureq} and set
\begin{equation}\label{definiofsigmai}
\sigma_i=\frac{1}{2\pi}\int_{\mathbb R^2}e^{u_i}dx\text{\, is finite},\quad i\in I.
\end{equation}
As did in \cite{CSW}, $u_i$ satisfies the asymptotic behavior of \eqref{new210}
where $m_i$ is defined by \eqref{new28}. By the fact $e^{u_i}\in L^1(\mathbb R^2)$ and
by using the Brezis-Merle type argument (see Lemma 4.1 in \cite{linzhang1}),
it can be proved that $u_i\in L^\infty_{loc}(\mathbb R^2)$ and then $u_i\in C^\infty(\mathbb R^2)$ by
further applications of standard elliptic regularity theorems.
Set $v_i$ to be
\begin{equation}
v_i(x)=\frac{1}{2\pi}\int_{\mathbb R^2}\log\frac{|x|}{|x-y|}e^{u_i(x)}dx.
\end{equation}
Then $u_i-\sum_{j\in I}a_{ij}v_j$ is a harmonic function in $\mathbb R^2$ and is bounded from above by $c\log|x|$ for large $x$.
Thus, $u_i(x)=\sum_{j\in I}a_{ij}v_j(x)+c_i$ for some constant $c_i$. Clearly, $(v_1,\cdots, v_n)$ satisfies \eqref{chipotshi}:
\begin{equation*}
-\Delta v_i=e^{u_i}=\mu_i e^{\sum_{j\in I} a_{ij}v_j}
\end{equation*}
for $\mu_i=e^{c_i}$. Therefore, $\sigma=(\sigma_1,\cdots,\sigma_n)$ satisfies \eqref{iffcon}.

\medskip

We note that in \cite{chanillo2} and \cite{CSW}, it has been shown that any component $u_i$ of a solution $u$ of \eqref{oureq} with $e^{u_i}\in L^1(\mathbb R^2)$
must be radially symmetric with respect to some point in $\mathbb R^2$, in paticular,
if $A$ is irreducible, then $u_i$, $i\in I$, are radially symmetric with respect to one common point.
If $A$ is not irreducible, then by \eqref{matrixa}, $u_i$, $i\in I_l$, is symmetric with respect to some $p_l$ for $l=1,2,\cdots,k$.
By replacing $u_i$ by $u_i(x+p_l)$, $i\in I_l$, we conclude that for any $\sigma$ satisfying \eqref{iffcon},
there exists a radial solution $u_i$ of \eqref{oureq} such that $\sigma_i=\frac{1}{2\pi}\int_{\mathbb R^2}e^{u_i}dx$.
In summary, Theorem B can be written as the following result for system \eqref{oureq}.

\medskip

\noindent
\emph{Theorem C: Let $A$ be the same as in Theorem B. Then a solution $u=(u_1,\cdots,u_n)$ of \eqref{oureq} exists with $\sigma_i=\frac{1}{2\pi}\int_{\mathbb R^2}e^{u_i}dx$
 if and only if $\sigma$ satisfies \eqref{iffcon}.
Furthermore, after a translation, the solution $u$ is radially symmetric with respect to the origin.}

\medskip

Set
\begin{equation}\label{energycon}
\mathcal{E}:=\{\sigma; \sigma_i>0, \, i\in I;\,\,  \Lambda_I(\sigma)=0;\,\, \Lambda_J(\sigma)>0 \,\, \forall\, \emptyset\subsetneq J\subsetneq I.\}.
\end{equation}
Let $u$ be a radial solution with $\sigma\in\mathcal{E}$, where $\sigma$ is given by \eqref{definiofsigmai}. Then it is clear to see that
\begin{equation*}
\tilde u_i(x)=u_i(\delta x)+2\log\delta,\quad i\in I,
\end{equation*}
is also a solution of \eqref{oureq} with the same $\sigma_i$. Thus, without loss of generality, we may assume $u$
satisfies $u_i(0)=\alpha_i$, $u_n(0)=0$, $1\leq i\leq n-1$, and let $\mathcal{B}\subseteq\mathbb R^{n-1}$ be the set of
initial values of solutions in \eqref{oureq}. Assume $({\bf H1})$ holds. In \cite{linzhang1}, the authors proved
\begin{equation}\label{515e1}
\mathcal{B} \mbox{ is open
and is homeomorphic to } \mathcal{E} \quad \mbox{(defined in (\ref{energycon}))}.
\end{equation}

In particular $\mathcal{E}$ is an open set in $\{\sigma;\,\, \Lambda_I(\sigma)=0\,\, \}$. If $a_{ii}>0$ for all $i\in I$, then it is not difficult
to see $\mathcal{B}=\mathbb R^{n-1}$. In general $\mathcal{B}$ might not be equal to $\mathbb R^{n-1}$, it is even not known whether $\mathcal{B}\neq \emptyset$.

In the following Lemma \ref{irredu} and Theorem \ref{calculusthm} we show that with ({\bf H2}), all submatrices of $A$ are irreducible, $\mathcal{B}\neq \emptyset$ and the condition $\Lambda_J(\sigma)>0$ is automatically satisfied if $\Lambda_I(\sigma)=0$.

\begin{lem}\label{irredu}
Let $A$ satisfy $({\bf H1})$ and $({\bf H2})$, then $a_{ij}>0$ and
$\max (a_{ii}, a_{jj})\le a_{ij}$ for all $i\neq j$ in $I$. In particular, all submatrices of $A$ are irreducible.
\end{lem}

\noindent{\bf Proof of Lemma \ref{irredu}:}

 \medskip

 \noindent{\bf Step one: If $a_{i_0i_0}>0$ for some $i_0$, then $a_{i_0j}\ge a_{i_0i_0}$ for all $j\neq i_0$.}

 \medskip

 Suppose $a_{11}>0$, let $\sigma_1=1,\sigma_2=..=\sigma_n=0$ and $m_i=\sum_{j=1}^na_{ij}\sigma_j$. Then $m_j=a_{1j}$ for all $j\in I$. Let $m=\min\{m_2,...,m_n\}$. We want to show that $m>0$. Indeed, if $m=0$, we let
 $$J=\{i\in I;\quad m_i=0. \}.$$
 Clearly $1\not \in J$, so for any $i\in J$, $\sigma_i=0$, which reads
 $$
 0=\sigma_i=\sum_{j\in J}a^{ij}m_j+\sum_{j\not \in J}a^{ij}m_j=\sum_{j\not \in J}a^{ij}m_j,\quad \forall i\in J.
 $$
Since $a^{ij}\ge 0$ for all $i\neq j$ (see $({\bf H2})$), we have $a^{ij}=0$ for all $i\in J$ and $j\not \in J$. After a permutation of the rows of
$A^{-1}$ (therefore the same permutation on the columns) $A^{-1}$ is of the form:
$$\left(\begin{array}{cc}
{\bf B_1} &  {\bf 0 }\\
{\bf 0}   &  {\bf B_2}
\end{array}
\right )
$$
which means $A$ is of the form
$$\left(\begin{array}{cc}
{\bf B_1^{-1}} &  {\bf 0 }\\
{\bf 0}   &  {\bf B_2^{-1}}
\end{array}
\right )
$$
after a permutation of its rows and columns. This is a contradiction to the irreducibility of $A$. Therefore we have proved $m>0$. Next we claim $m\ge m_1$. If this is not true, we have $m<m_1$ and we let
$$J:=\{i\in I;\quad m_i=m\quad \}. $$
Then by our assumption, $J\neq \emptyset$ and $1\not \in J$. Thus, for any $i\in J$, the fact $\sigma_i=0$ yields
\begin{eqnarray}\label{521e1}
0&=&\sigma_i=\sum_{j\in J}a^{ij}m_j+\sum_{j\not \in J}a^{ij}m_j\\
&\ge &m\sum_{j\in J}a^{ij}+m\sum_{j\not \in J}a^{ij}\nonumber \\
&=&m\sum_{j\in I}a^{ij}\ge 0,\quad \forall i\in J.\nonumber
\end{eqnarray}
Note that we have used $(\bf {H2})$ in both inequalities. We see that the second line of (\ref{521e1}) is a strict inequality unless
$a^{ij}=0$ for all $i\in J$ and $j\not \in J$. Thus, (\ref{521e1}) yields $a^{ij}=0$ for $i\in J$ and $j\not\in J$, a contradiction to the irreducibility of $A$. Step one is established.

\medskip

\noindent{\bf Step two: $a_{ij}>0$ for all $i\neq j$}

\medskip

We prove step two by contradiction. Suppose there exist $i\neq j$ such that $a_{ij}=0$. By step one, we have $a_{ii}=a_{jj}=0$. Without loss of generality we assume $i=1$ and $j=2$. We can also assume $a_{13}>0$, because the invertibility of $A$ implies $a_{1i}>0$ for some $i\ge 3$. We can apply a permutation on $A$ to move the positive entry to the third row. Thus the matrix $A$ is of the following form after a permutation of rows and columns:
$$\left(\begin{array}{cccc}
0 & 0 & a_{13} & \ldots \\
0 & 0 & a_{23} & \ldots \\
a_{13} & a_{23} & a_{33} & \ldots\\
\vdots & \vdots & \vdots & \ddots
\end{array}
\right )
$$
Let $\sigma_1=1, \sigma_2=0, \sigma_3=1$ and $\sigma_4=..=\sigma_n=0$ and $m_i=\sum_{j=1}^na_{ij}\sigma_j$. Clearly $m_1=a_{13}$, $m_3\ge a_{13}$. Let $m=\min\{m_2,m_4,..,m_n\}$. We first claim $m>0$. Suppose this is not the case, let $J=\{i\in I;\quad m_i=0\}$. Obviously $1,3\not \in J$, which implies $\sigma_j=0$ for all $j\in J$. Using the same argument as in step one, we have $a^{ij}=0$ for all $i\in J$ and all $j\not \in J$, thus a contradiction to the irreducibility of $A$.

Next we claim $m\ge m_1$. Suppose this is not true, then $m<m_1$. Observe that $m_3\ge m_1$ so $m<m_3$. Let
$$J=\{i\in I;\quad m_i=m\}. $$
Then $1,3\not \in J$, so $\sigma_i=0$ for all $i\in J$. Then \eqref{521e1} yields
\begin{equation*}
0=\sigma_i\geq m\sum_{j\in I}a^{ij}+2(m_1-m)a^{13}>0,
\end{equation*}
which is a contradiction. Therefore $m\ge m_1$ is proved. In particular, $m_2=a_{23}\ge m_1$, which gives $a_{23}\ge a_{13}$. Since we can switch the first two rows of $A$ the same argument gives $a_{13}\ge a_{23}$, consequently $a_{13}=a_{23}$. We use $a_{13}$ to represent any nonzero entry on the first row, so we have proved that the first two rows are identical, a contradiction to the invertibility of $A$. Lemma \ref{irredu} is established. $\Box$

\begin{thm}\label{calculusthm}
Let $A$ satisfy $(\bf H1)$ and $(\bf H2)$. Suppose $\sigma=(\sigma_1,..,\sigma_n)$ has positive components and $\Lambda_I(\sigma)=0$. Then
$\Lambda_J(\sigma)>0$ for all $\emptyset \subsetneq J \subsetneq I$.
\end{thm}

\noindent{\bf Proof of Theorem \ref{calculusthm}:}

First we assume $\mathcal{E}\neq \emptyset$ (which will be proved in Lemma \ref{enotempty}) and prove
\begin{eqnarray}\label{bryE}
\mbox{If }\tilde \sigma\in \partial \mathcal{E}, \mbox{ then } \exists J\subsetneq I,\quad  \mbox{ such that } \tilde \sigma_i>0 \quad
 \forall i\in J; \\
 \sigma_i=0 \,\, \forall i\not \in J;\quad \mbox{and $\tilde \sigma$ satisfies }
\,\, \Lambda_J(\tilde \sigma)=0. \nonumber
\end{eqnarray}
\noindent{\bf Proof of (\ref{bryE}):}

Let $\sigma^k$ be a sequence of points in $\mathcal{E}$ that tends to $\tilde \sigma\in \partial \mathcal{E}$. Let $u^k=(u^k_1,..,u^k_n)$ be global, radial solutions that correspond to $\sigma^k$. Without loss of generality we assume $u_1^k(0)=\max_{i\in I}u^k_i(0)=0$. Since $e^{u_i^k}\le 1$, by the standard elliptic estimates, there exists a subsequence of $u^k$ (still denoted by $u^k$) such that parts of $u^k$ converge in $C^{\infty}_{loc}(\mathbb R^2)$. There are two cases of the convergence of $u^k$ to be discussed separately.

\medskip

{\bf Case one:}
$u^k$ converges to a global solution $v=(v_1,..,v_n)$ which satisfies the system (\ref{oureq}) of $n$ equations.

\medskip

We claim this case can not happen.
Indeed, let
$$\sigma_i^k:=\frac 1{2\pi}\int_{\mathbb R^2}e^{u_i^k} \mbox{ and } \sigma_i:=\frac 1{2\pi}\int_{\mathbb R^2}e^{v_i}, \quad i\in I. $$
Clearly
\begin{equation}\label{521e6}
\tilde \sigma_i:=\lim_{k\to \infty}\sigma_i^k\ge \sigma_{i}.
\end{equation}
Note that $\tilde \sigma=(\tilde \sigma_1,..,\tilde \sigma_n)$ and
$\sigma:=(\sigma_{1},..,\sigma_{n})$ both satisfy $\Lambda_I=0$ and we also have $\sum_{j\in I}a_{ij}\sigma_{j}>2$ for $i\in I$ because $\int_{\mathbb R^2}e^{v_i}<\infty$. Taking the difference on the two equations $\Lambda_I(\tilde \sigma)=0$ and $\Lambda_I(\sigma)=0$, we arrive at
\begin{equation}\label{521e7}
\sum_{j\in I}(\sum_{i\in I}a_{ij}\tilde \sigma_i-2)(\tilde \sigma_j-\sigma_j)+
\sum_{j\in I}(\sum_{i\in I}a_{ij}\sigma_i-2)(\tilde \sigma_j-\sigma_j)=0.
\end{equation}
For each $i\in I$, $\sum_{j\in I}a_{ij}\tilde \sigma_j\ge \sum_{j\in I}a_{ij}\sigma_j>2$. Combining this fact with (\ref{521e6}) and (\ref{521e7}) we have
$\tilde \sigma_i=\sigma_{i}$ for all $i\in I$.
Thus $\tilde \sigma\in \mathcal{E}$, which is an open subset of the hypersurface $\Lambda_I=0$ (see (\ref{515e1})), a contradiction to the assumption that $\tilde \sigma \in \partial \mathcal{E}$.

\medskip

\noindent{\bf Case two:} There exists $K\subsetneq I$ such that $u_j^k$ converges for $j\in K$ and $u_j^k(r)\to -\infty$ uniformly in any bounded set of $[0,\infty)$ for $j\not \in K$.

\medskip

Let $l=|K|$, clearly $u_1^k$ must converge, so without loss of generality, we assume that the first $l$ components of $u^k$ converge to $v=(v_1,..,v_l)$ which satisfies

\begin{equation}\label{518e15}
-\Delta v_i=\sum_{j=1}^la_{ij}e^{v_j},\quad i=1,..,l, \quad \mbox{ in } \mathbb R^2
\end{equation}
and it is easy to show
$$\tilde \sigma_i:=\lim_{k\to \infty}\sigma_i^k\ge \sigma_i:=\frac 1{2\pi}\int_{\mathbb R^2}e^{v_i},\quad i=1,..,l.$$
Since $\int_{\mathbb R^2}e^{v_i}<\infty$, by \eqref{new27}, we have
\begin{equation}\label{514e1}
\sum_{j=1}^la_{ij}\sigma_{j}>2,\quad i=1,..,l.
\end{equation}
Let $\sigma:=(\sigma_{1},..,\sigma_{l},0,..,0)$ and $\tilde \sigma=(\tilde \sigma_1,..,\tilde \sigma_n)$ be the limit of $\sigma^k$. We have proved
that $\tilde \sigma_i\ge \sigma_{i}$ for $1\le i\le l$. Although $(a_{ij})_{l\times l}$ may not be invertible, we still have the Pohozaev identity:
\begin{equation}\label{partialpi}
\sum_{i,j=1}^la_{ij}\sigma_{i}\sigma_{j}=4\sum_{i=1}^l\sigma_{i}.
\end{equation}
See Theorem C. Now we claim
\begin{equation}\label{515e2}
\tilde \sigma_i=\sigma_{i},\quad i=1,..,l;\quad \tilde \sigma_i=0\quad \mbox{ for }\quad i>l.
\end{equation}

Clearly (\ref{bryE}) follows from (\ref{515e2}) and (\ref{partialpi}).
The proof of (\ref{515e2}) relies heavily on $({\bf H2})$.
Set
$$m_{i}=\sum_{j=1}^na_{ij}\sigma_{j},\quad i\in I. $$
We want to prove
\begin{equation}\label{224thm21}
m_{i}>2\quad\text{for all \,}i\in I.
\end{equation}
Obviously this is true for $i=1,..,l$, so we need to prove $m_{i}>2$ for
$i>l$. To see this, we first observe that $m_{i}>0$ for all $i$ because $a_{ij}>0$, $i\not=j$, by Lemma \ref{irredu}.

The proof of $m_{i}>2$ can be obtained by the similar argument of Lemma \ref{irredu}. Let
$$m=\min\{m_{l+1},m_{l+2},..,m_{n}\}.$$
Suppose $m\leq 2$, then let $J=\{i\in I;\, m_i=m\}$ and $m^*=\min_{i\in J^c}\{m_i\}$.
Same as before, $1,..,l\not \in J$, and $\sigma_i=0$ for all $i\in J$. So, \eqref{521e1} yields
\begin{equation*}
0=\sigma_i\geq m\sum_{j\in I}a^{ij}+(m^*-m)\sum_{j\not\in J^c}a^{ij}\geq 0,
\end{equation*}
and then $a^{ij}=0$ for $i\in J$ and $j\not\in J^c$, which yields a contradiction to the irreducibility of $A$. Hence, $m_i>2$ is proved.

\medskip

Now we finish the proof of (\ref{515e2}). Certainly, (\ref{partialpi}) can be written as $\Lambda_I(\sigma)=0$ (the last $n-l$ components of $\sigma$ are $0$). $\tilde \sigma$ also satisfies $\Lambda(\tilde \sigma)=0$. Besides we have $\tilde \sigma_i\ge \sigma_i$ and $\tilde m_i, m_i>2$ for all $i\in I$.
Using \eqref{521e7}, we obtain $\tilde \sigma_i=\sigma_i$ for all $i$.
(\ref{515e2}) and (\ref{bryE}) are established.

\medskip

Finally, for any $\sigma$ on the surface $\Lambda_I=0$ with all the components positive, we claim that $\sigma\in \mathcal{E}$. Suppose $\sigma\not \in \mathcal{E}$.
Let $\sigma_E\in \mathcal{E}$, since $\Lambda_I(\sigma)=0$ is connected, we can find a path $\Gamma(t) (0\le t\le 1)$ on $\Lambda_I=0$ that connects $\sigma$ and $\sigma_E$ ($\Gamma(0)=\sigma_E,\Gamma(1)=\sigma$). Here we require all the components of $\Gamma(t)$ ($0\le t\le 1$) be positive. Because $\sigma\not \in \mathcal{E}$, there exists $t_0\in [0,1]$ such that $\Gamma(t_0)\in \partial \mathcal{E}$. But it yields a contradiction to (\ref{bryE}) because no component of $\Gamma(t_0)$ is zero.
Theorem \ref{calculusthm} is established. $\Box$

\bigskip

In the proof of Theorem \ref{calculusthm} we assumed $\mathcal{E}\neq \emptyset$, which is established in the following lemma.

\begin{lem}\label{enotempty}
$\mathcal{E}$ is not empty.
\end{lem}

\noindent{\bf Proof of Lemma \ref{enotempty}:}

Let $\xi_i=\sum_{j=1}^na^{ij}$ for $i\in I$. From ${\bf (H2)}$ we know $\xi_i\ge 0$ for all $i\in I$. First consider the case that all $\xi_i$ are positive:

\medskip

\noindent{\bf Case one: $\xi_i>0$ for all $i\in I$}

\medskip

Using the properties of inverse matrices
\begin{equation}\label{518e8}
\sum_{j=1}^na_{ij}\xi_j=\sum_{j=1}^na_{ij}\sum_{k=1}^na^{jk}=\sum_{k=1}^n\delta_{ik}=1.
\end{equation}
Let $\sigma_i=4\xi_i$ for $i\in I$. Then by the assumption in this case $\sigma_i>0$ for all $i\in I$ and direct computation shows
\begin{eqnarray}\label{518e9}
&&\Lambda_I(\sigma)=4\sum_{i\in I}\sigma_i-\sum_{i,j\in I}a_{ij}\sigma_i\sigma_j\\
&=&16(\sum_{i\in I}\xi_i-\sum_{i,j\in I}a_{ij}\xi_i\xi_j)
=16\sum_{i\in I}\xi_i(1-\sum_{j=1}^na_{ij}\xi_j)=0.\nonumber
\end{eqnarray}
For any nonempty $J\subsetneq I$, without loss of generality $J={1,..,l}$ for some $l<n$, easy to see
\begin{eqnarray*}
\Lambda_J&=&4\sum_{i=1}^l\sigma_i-\sum_{i,j=1}^la_{ij}\sigma_i\sigma_j\\
&=&16(\sum_{i=1}^l\xi_i-\sum_{i,j=1}^la_{ij}\xi_i\xi_j)\\
&=&16\sum_{i=1}^l\xi_i(1-\sum_{j=1}^la_{ij}\xi_j)=16\sum_{i=1}^l(\sum_{j=l+1}^na_{ij}\xi_j)\xi_i.
\end{eqnarray*}
Clearly $\Lambda_J\ge 0$. Since by Lemma \ref{irredu} $a_{ij}>0$ for $1\le i\le l$ and $l+1\le j\le n$ we have
$\Lambda_J>0$. Therefore $\sigma\in \mathcal{E}$ and $\mathcal{E}\neq \emptyset$.

\medskip

\noindent{\bf Case Two: There exists $\xi_i=0$}

\medskip

First we observe that it is not possible to have all $\xi_i=0$ because otherwise adding all the rows of $A^{-1}$ to the first row would make all the entries of the first row $0$, a contradiction to the invertibility of $A^{-1}$.
Without loss of generality we assume $\xi_{l+1}=...=\xi_n=0$ and $\xi_1,..,\xi_l>0$.
Let $J=\{1,..,l\}$ then we claim
\begin{equation}\label{518e1}
 \Lambda_J(\sigma)=0\quad \mbox{and}\quad \Lambda_{J_1}>0 \quad \forall \emptyset\subsetneq J_1\subsetneq J.
 \end{equation}
From (\ref{518e9}) we see easily $\Lambda_J=0$. For $\emptyset\subsetneq J_1\subsetneq J$, without loss of generality we assume $J_1=\{1,..,l_1\}$ with $l_1<l$. Similar to case one (using $\sum_{j=1}^la_{ij}\xi_j=1$ and $\sigma_i=4\xi_i$ for all $i=1,..,l$)
$$\Lambda_{J_1}=16\sum_{i=1}^{l_1}(\sum_{j=l_1+1}^la_{ij}\xi_j)\xi_i. $$
Thus $\Lambda_{J_1}>0$ is an immediate consequence of Lemma \ref{irredu}. (\ref{518e1}) is established.

\medskip

Let $\tilde A=(a_{ij})_{i,j\in J}$. Although $\tilde A$ may not be invertible, Theorem B can still be applied to conclude that there exists a radially symmetric solution $u(r)=(u_1(r),..,u_l(r))$ to
\begin{equation}\label{518e4}
\left\{\begin{array}{ll}
\Delta u_i+\sum_{j=1}^{l}a_{ij}e^{u_j}=0,\quad \mathbb R^2,\quad i=1,..,l,\\
\\
\frac{1}{2\pi }\int_{\mathbb R^2}e^{u_i}=\sigma_i=4\xi_i,\quad i=1,..,l.
\end{array}
\right.
\end{equation}
Since $e^{u_i}\in L^1(\mathbb R^2)$, we have
\begin{equation}\label{518e6}
\sum_{j=1}^{l}a_{ij}\sigma_j>2, \quad i=1,..,l.
\end{equation}
It has been discussed in the proof of Theorem \ref{calculusthm} that $({\bf H2})$ implies
\begin{equation}\label{518e30}
\sum_{j=1}^{l}a_{ij}\sigma_j>2,\quad \forall i\in I.
\end{equation}
which leads to
\begin{equation}\label{518e10}
\frac 1{2\pi}\int_{B_R}\sum_{j=1}^{l}a_{ij}e^{u_j}>2+\delta,\quad i\in I
\end{equation}
for some $\delta>0$ and $R>1$.
Now we construct a sequence of functions $u^{\epsilon}=(u^{\epsilon}_1,..,u^{\epsilon}_n)$ as follows:
$$\left\{\begin{array}{ll}
(u_i^{\epsilon})''(r)+\frac 1r(u_i^{\epsilon})'(r)+\sum_{j=1}^na_{ij}e^{u_j^{\epsilon}(r)}=0,\quad 0<r<\infty,\quad i\in I,\\
\\
u_i^{\epsilon}(0)=u_i(0),\quad i=1,...,l,\\
\\
u_i^{\epsilon}(0)=\log \epsilon,\quad i=l+1,..,n.
\end{array}
\right.
$$
By standard ODE existence theory, solution $u^{\epsilon}$ is well defined for all $r>0$. Since $u^{\epsilon}_i(r)$ are decreasing functions, it is easy to see that as $\epsilon$ tends to $0$, $u^{\epsilon}_i$ ($i>l$) tends to $-\infty$ over $[0, R]$ for any fixed $R>0$. Therefore, the first $l$ components of $u^{\epsilon}$ converge uniformly to the corresponding components of $u$ over any fixed $[0, R]$. In particular for the $R$ in
(\ref{518e10}) we have, for $\epsilon$ sufficiently small
\begin{equation}\label{518e33}
\frac{1}{2\pi}\int_{B_R}\sum_{j=1}^na_{ij}e^{u_j^{\epsilon}}\ge
\frac 1{2\pi}\int_{B_R}\sum_{j=1}^{l}a_{ij}e^{u^{\epsilon}_j}>2+\delta/2,\quad \forall i\in I,
\end{equation}
which implies
$$(u_i^{\epsilon})'(r)=-\frac{1}r\int_0^{r}\sum_{j=1}^na_{ij}e^{u_j^{\epsilon}(s)}sds<-\frac{2+\delta/2}{r},\quad r>R. $$
Thus
$$u_i^{\epsilon}(r)\le -(2+\delta/2)\log r+O(1),\quad \mbox{ for } r>R $$
where $O(1)$ is a constant independent of $r$ but may depend on $\epsilon$. Hence
for $\epsilon>0$ small, $\int_{\mathbb R^2}e^{u_i^{\epsilon}}<\infty$ for all $i\in I$. Thus $\mathcal{E}\neq \emptyset$. Case two and Lemma \ref{enotempty} are established. $\Box$

\begin{rem}\label{519rem1}
For the $2\times 2$ case, let $A=\left(\begin{array}{cc}
a & c\\
c & b
\end{array}
\right )
$ be an non-negative and irreducible matrix (which means $c>0$). It
can be proved that the conclusion of Theorem \ref{calculusthm} holds for $A$ if and only if $\max(a,b)\le 2c$. However, $A$ satisfies ({\bf H1}) and ({\bf H2}) if and only if $\max(a,b)\le c$ and $c^2\neq ab$.
\end{rem}

\section{The case of partial blowup}

In this section we consider the following case: Let $u^k=\{u_1^k,..,u_n^k\}$ be a sequence of solutions of
\begin{equation}\label{uikeq}
-\Delta u_i^k=\sum_{j=1}^na_{ij}h_j^k(x)e^{u_j^k},\quad B_1,\quad i\in I
\end{equation}
where $h^k=(h_1^k,..,h_n^k)$ is a sequence of positive functions with uniformly bounded $C^1$ norm:
\begin{equation}\label{hika}
\frac 1C\le h_i^k(x)\le C, \quad |\nabla h_i^k|(x)\le C,\quad \forall x\in B_1,\quad \forall i\in I.
\end{equation}
Suppose that $u^k$ blows up only at $0$:
\begin{equation}\label{219e1}
M_k=\max_{i\in I, x\in B_1}u_i^k(x)\to \infty,\quad
\max_{i\in I} u_i^k\le C(K)\,\, \forall K\subset\subset B_1\setminus \{0\}.
\end{equation}
Let
\begin{equation}\label{sigmdef}
\sigma_i^k=\frac 1{2\pi}\int_{B_1}h_i^ke^{u_i^k};
\quad m_i^k=\sum_{j=1}^na_{ij}\sigma_j^k,\quad i\in I.
\end{equation}
Without loss of generality, we assume $u_1^k(0)=M_k$, and set
\begin{equation}\label{videf}
v_i^k(y)=u_i^k(\delta_k y)+2\log \delta_k,\quad\text{where \,}\delta_k=e^{-\frac{1}{2}u_1^k(0)}.
\end{equation}
By elliptic estimates, we can show that there is $J\subseteq I$ such that $\{v_i^k\}_{i\in J}$ converges, in $C^2_{loc}(\mathbb R^2)$, to $\{v_i\}_{i\in J}$.
We may assume $J=\{1,2,\cdots, l\}$. Then $v=(v_1,\cdots,v_l)$ solves the following subsystem:
\begin{equation}\label{equofvi}
-\Delta v_i=\sum_{j=1}^l a_{ij}\mathcal{H}_j e^{v_j},\quad i=1,\cdots, l, \quad\text{in \,}\mathbb R^2,
\end{equation}
where $\mathcal{H}_i=\lim_{k\to \infty}h_i^k(0)$.
In this case, we have $v_i^k(x)\rightarrow-\infty$ for $i\not \in J$ as $k\rightarrow+\infty$ in any compact set of $\mathbb R^2$.

\medskip

We also assume $u^k$ to have bounded oscillation on $\partial B_1:$
\begin{equation}\label{219e2}
|u_i^k(x)-u_i^k(y)|\le C,\quad \forall x, y\in \partial B_1,\quad i\in I
\end{equation} and uniformly bounded energy in $B_1$.
\begin{equation}\label{219e4}
\int_{B_1}h_i^ke^{u_i^k}\le C, \quad i\in I.
\end{equation}

\medskip

When $l=n$, it was proved by the authors \cite{linzhang1} that there exists an entire solution $U^k=(U_1^k,\cdots, U_n^k)$ of
\begin{equation*}
\Delta U_i^k+\sum{a_{ij}}\mathcal{H}_j e^{U_i^k}=0\quad\text{in \,}\mathbb R^2
\end{equation*}
such that
\begin{equation*}
|u_i^k(x)-U_i^k(x)|\leq C\quad\text{for \,}|x|\leq\frac{1}{2}\quad\text{and}\quad i\in I.
\end{equation*}
An immediate consequence is the following estimate:
\begin{equation}\label{estimateofu}
u_i^k(x)=-\frac{m_i^k-2}{2}M_k+O(1)\quad\text{for \,}|x|=\frac{1}{2}.
\end{equation}

\medskip

In this section, we want to extend the estimate \eqref{estimateofu} to the case when
there are only $l$ components of $v^k$, $l<n$, which converge to an entire solution of \eqref{equofvi}.

\begin{prop}\label{blowuplocal}
Let $h^k=(h_1^k,\cdots, h_n^k)$ satisfy (\ref{hika}), and let
$u^k$ be solutions of (\ref{uikeq}) such that (\ref{219e1}), (\ref{equofvi}), (\ref{219e2}) and (\ref{219e4}) hold.
Let
$\sigma_i=\lim_{k\rightarrow+\infty}\sigma_i^k$ and $m_i=\lim_{k\rightarrow+\infty}m_i^k$,
where $\sigma_i^k$ and $m_i^k$ are given by \eqref{sigmdef}.
Then
\begin{enumerate}
\item $4\sum_{i=1}^l\sigma_i=\sum_{i,j=1}^la_{ij}\sigma_i\sigma_j. $
\item $\sigma_1,..,\sigma_l>0$, $\sigma_{l+1}=...=\sigma_n=0$.
\item $m_i>2 \quad \forall i\in I$.
\item For $1\le i\le l$
\begin{equation}\label{511e1}
|u_i^k(x)+\frac{m_i^k-2}2M_k|\le C(\epsilon),\quad \forall x\in B_1\setminus \bar B_{\epsilon},
\end{equation}
\item For $l+1\le i\le n$
\begin{equation}\label{219e7}
|u_i^k(x)+\frac{m_i^k-2}2M_k+(M_k-u_i^k(0))|\le C(\epsilon),\quad \forall x\in B_1\setminus \bar B_{\epsilon}.
\end{equation}
\end{enumerate}
\end{prop}

\noindent{\bf Proof of Proposition \ref{blowuplocal}:}

\medskip

This proof is divided into a few steps.

\medskip

\noindent{\bf Step one: $m_i>2$ for all $i\in I$.}

\medskip

Let
$$\sigma_{iv}=\frac 1{2\pi}\int_{\mathbb R^2}\mathcal{H}_i(0)e^{v_i},\quad i=1,2,..,l. $$
Then Theorem C and \eqref{new27} imply
\begin{equation}\label{piv1}
\left\{\begin{array}{ll}\sum_{i,j=1}^la_{ij}\sigma_{iv}\sigma_{jv}=4\sum_{i=1}^l\sigma_{iv}\\
\\
\sum_{j=1}^la_{ij}\sigma_{jv}>2,\quad i=1,...,l.
\end{array}
\right.
\end{equation}
Obviously, $\sigma_i\ge \sigma_{iv}$ for $i=1,..,l$. Let $\sigma_{iv}=0$ for $i=l+1,..,n$ and
$$m_{iv}=\sum_{j=1}^na_{ij}\sigma_{jv},\quad i=1,..,n. $$
We claim that $m_{iv}>2$ for all $i=1,..,n$. For $i\leq l$, this is already known by \eqref{piv1}. The proof of $m_{iv}>2$ for
$i>l$ is the same as \eqref{224thm21} in the proof of Theorem \ref{calculusthm}.
Since $\sigma_i\geq \sigma_{iv}$, we also have $m_i>2$.

\medskip
\noindent{\bf Step two: $u_i^k(x)\to -\infty$ over any compact subset of $B_1\setminus \{0\}$ for all $i=1,2,..,n$.}

\medskip
We first show that $u_i^k$ tends to $-\infty$ on $\partial B_1$.
Let
$G(x,y)$ be the Green's function with respect to the Dirichlet condition on $\partial B_1$, the Green's representation formula
gives
$$u_i^k(x)\ge \int_{B_1}G(x,\eta)\bigg (\sum_{j=1}^la_{ij}h_j^k(\eta)e^{u_j^k(\eta)}\bigg )d\eta+\min_{\partial B_1}u_i^k. $$
If $\min_{\partial B_1}u_i^k\ge -C$ for some $C>0$, we use
$$\sum_{j=1}^la_{ij}h_j^ke^{u_j^k}\to 2\pi m_i\delta_0\quad \mbox{ in measure} $$ where $\delta_0$
is the Dirac mass at $0$, and
$$G(x,\eta)\ge -\frac{1}{2\pi}\log |x-\eta |-C_1\quad\text{for \,}|x|\leq \frac{1}{2} $$
to obtain
$$\lim_{k\to \infty}e^{u_i^k(x)}\ge C_2|x|^{-2-\epsilon_1}\quad\text{for \,}0<|x|\leq\frac{1}{2}\quad\text{and for some \,}\epsilon_1>0, $$
where $C_1$ and $C_2$ are two positive constants. This is a contradiction to $\int_{B_1}e^{u_i^k}\le C$. Therefore we have proved $\min_{\partial B_1}u_i^k\to
-\infty$. Since $u_i^k$ is bounded from above over any compact subset of $B_1\setminus \{0\}$, standard elliptic estimate implies
that $u_i^k\to -\infty$ uniformly in any compact subset of $B_1\setminus \{0\}$.

\medskip

\noindent{\bf Step three: }$\sigma_{l+1}=...=\sigma_n=0$.

\medskip

Let $w_i^k=\sum_{j\in J}a^{ij}u_j$. Then $w_i^k$ satisfies
\begin{equation*}
-\Delta w_i^k=h_i^k e^{\sum_{j\in I}a_{ij}w_j^k}.
\end{equation*}
We can apply the Pohozaev identity to the above equation, i.e.,
multiplying $x\cdot\nabla (\sum a_{ij}w_j^k)$ and applying integration by parts, we obtain by passing $k\rightarrow+\infty$,
$$\sum_{i,j=1}^na_{ij}\sigma_i\sigma_j=4\sum_{i=1}^n\sigma_i. $$
On the other hand, $\sigma_v=(\sigma_{1v},...,\sigma_{lv},0,..,0)$ also satisfies the equality above. Besides, we also have
$$\sigma_i\ge \sigma_{iv},\quad i=1,..,l,\quad \sigma_i\ge 0,\quad i=l+1,..,n. $$
By taking the difference of the above equations satisfied by $\sigma$ and $\sigma_v$, we arrive at the same as \eqref{521e7}:
\begin{equation*}
\sum_{j\in I}\left(\sum_{i\in I}a_{ij}\sigma_i-2\right)(\sigma_j-\sigma_{jv})+
\sum_{j\in I}\sum_{i\in I}(a_{ij}\sigma_{iv}-2)(\sigma_j-\sigma_{jv})=0.
\end{equation*}
Recall $m_i=\sum a_{ij}\sigma_j>2$ and $m_{iv}=\sum a_{ij}\sigma_{jv}>2$. Then the above identity leads to $\sigma=\sigma_v$, i.e.,
$\sigma_i=\sigma_{iv}$ for $i\leq l$ and $\sigma_i=0$ for $i>l$.

\medskip

\noindent{\bf Step four: $u_i^k(x)+2\log |x|\le C$}

\medskip

By scaling, it is equivalent to proving
\begin{equation}\label{223e2}
v_i^k(y)+2\log |y|\le C, \quad |y|\le \delta_k^{-1}.
\end{equation}
Since only $\{v_1^k,..,v_l^k\}$ converges to $\{v_1,..,v_l\}$ in $C^2_{loc}(\mathbb R^2)$, it implies that $v_{l+1}^k,...,v_n^k$ tend to
$-\infty $ in any compact subset of $\mathbb R^2$.
Suppose (\ref{223e2}) does not hold, there exists $y_k\to \infty$ such that, along a subsequence
\begin{equation}\label{522e5}
v_{i_0}^k(y_k)+2\log |y_k|=\max_{|y|\le \delta_k^{-1},i\in I}(v_i^k(y)+2\log |y|)\to +\infty.
\end{equation}
It is easy to see that $|y_k|\delta_k\to 0$ as $u_i^k$ is bounded above in any compact subset of $B_1\setminus \{0\}$.
Consider $y\in B(y_k,|y_k|/2)$, from
$$v_i^k(y)+2\log |y|\le v_{i_0}^k(y_k)+2\log |y_k|,\quad i\in I $$
 we have
\begin{equation}\label{223e3}
v_i^k(y)\le v_{i_0}^k(y_k)+2\log 2\quad\text{for \,}y\in B(y_k,|y_k|/2)\quad\text{and}\quad i\in I.
\end{equation}
On one hand, the fact $|y_k|\to \infty$ implies
\begin{equation}\label{522e6}
\lim_{k\to \infty} \frac{1}{2\pi}\int_{B(0, |y_k|/4)}h_i^k(\delta_k\cdot)e^{v_i^k}\ge \sigma_{iv},\quad i=1,2,..,l.
\end{equation}
On the other hand, we can set
$$w_i^k(y)=v_i^k(y_k+r_ky)-v_{i_0}^k(y_k),\quad \mbox{ for } |y|<R_k, \quad i\in I, $$
where $r_k=e^{-\frac 12v_{i_0}^k(y_k)}$ and $R_k=r_k|y_k|/2\to \infty$ as $k\to \infty$ by (\ref{522e5}). By (\ref{223e3}), $w_i^k(y)$ is bounded from above by $2\log 2$, and $w_{i_0}^k(0)=0$.
From the standard estimate of elliptic equations, there is $J\subseteq I$ such that $w_i^k(y)$ converges in $C^2_{loc}(\mathbb R^2)$ for $i\in J$ and
$w_i^k(x)\rightarrow -\infty$ for $i\not\in J$ in any compact set of $\mathbb R^2$. Clearly, $i_0\in J$. In particular,
$$\int_{|y|<R_k}e^{w_{i_0}^k}\ge \delta>0 \mbox{ for some }\delta>0. $$
Thus, together with (\ref{522e6}) we have
$$\sigma_{i_0}=\int_{B_1}h_{i_0}^ke^{u_{i_0}^k}\ge \sigma_{i_0v}+\delta>\sigma_{i_0v},$$
a contradiction to step three. Hence, step four is established.

\medskip

\noindent{\bf Step five: Estimate of the decay of $v_i^k$}.

\medskip

First we choose $R>>1$ such that
\begin{equation}\label{512e1}
\frac 1{2\pi}\int_{B_{R/2}} \sum_{j=1}^la_{ij}h_j^k(\delta_ky)e^{v_j^k(y)}dy>2+\delta, \quad i\in I.
\end{equation}
for some $\delta>0$. Note that (\ref{512e1}) obviously holds for $1\le i\le l$ because of of the convergence from $(v_1^k,..v_l^k)$ to $(v_1,..,v_l)$. Then by step one, it also holds for $i>l$.
In this step we study the behavior of $v_i^k$ for $|y|>R$ ($i\in I$).  For each $R<|y|<\delta_k^{-1}/2$, let
$$\hat v_i^k(z)=v_i^k(|y|z)+2\log |y|, \quad \frac 12<|z|<2. $$
By step four $\hat v_i^k\le C_1$ for some $C_1$ over $B_2\setminus B_{1/2}$. Consider the equation for $\hat v_i^k$ in $B_2\setminus B_{1/2}$:
$$-\Delta \hat v_i^k(z)=\sum_{j=1}^n a_{ij}h_j^k(\delta_k|y|z)e^{\hat v_j^k} \quad B_2\setminus B_{\frac 12}. $$
Let $f_i^k$ satisfy
$$\left\{\begin{array}{ll}
-\Delta f_i^k(z)=\sum_{j=1}^n a_{ij}h_j^k(\delta_k|y|z)e^{\hat v_j^k} \quad B_2\setminus B_{1/2}, \\
\\
f_i^k(z)=0,\quad \mbox{on}\quad \partial B_{1/2}\cup \partial B_2.
\end{array}
\right.
$$
Clearly $f_i^k\ge 0$ in $B_2\setminus B_{1/2}$. As a result of the upper bound of $\hat v_i^k$, we have
\begin{equation}\label{223e4}
0\le f_i^k(z)\le C_2, \quad z\in B_2\setminus B_{1/2}.
\end{equation}
Obviously, $\hat v_i^k-C_1-f_i^k$ is a non-positive harmonic function. Hence the Harnack inequality holds:
$$-\min_{\partial B_1}\left(\hat v_i^k-C_1-f_i^k\right)\le C\left(-\max_{\partial B_1}\left(\hat v_i^k-C_1-f_i^k\right)\right). $$
Equivalently,
$$\max_{\partial B_1}\hat v_i^k\le \frac 1{C}\min_{\partial B_1}\hat v_i^k+C_2. $$
Going back to $v_i^k$, we have
\begin{equation}\label{223e5}
\max_{\partial B_r}v_i^k(y)\le \frac{1}C\min_{\partial B_r}v_i^k+(-2+\frac 2C)\log r+C_2,
\end{equation}
for $R<r\le \delta_k^{-1}$ and $i=1,..,n$.

Let $\bar v_i^k(r)=\frac 1{2\pi r}\int_{\partial B_r}v_i^k$ be the average of $v_i^k$ on $\partial B_r$.
$$(\bar v_i^k)'(r)=\frac{1}{2\pi r}\int_{B_r}\Delta v_i^k. $$
By  (\ref{512e1}) we have
$$(\bar v_i^k)'(r)<-\frac{2+\delta}r,\quad r>R$$
for some $\delta>0$. So for $r>R$,
$$\bar v_i^k(r)\le -(2+\delta)\log r+\bar v_i^k(R). $$
By (\ref{223e5}), we have, for $i\in I$
\begin{eqnarray}\label{223e6}
v_i^k(y)&\le &\frac 1C\bar v_i^k(r)+(-2+\frac 2C)\log |y| +C_2,\\
&\le & -(2+\frac{\delta}C)\log |y|+C_2+\bar v_i^k(R). \nonumber
\end{eqnarray}
Note that for $i>l$, $\bar v_i^k(R)\to -\infty$. Thus
\begin{equation}\label{223e7}
v_i^k(y)\le -R_k-(2+\frac{\delta}C)\log |y|,\quad \forall i>l
\end{equation}
for some $R_k\to \infty$ as $k\to \infty$.

\medskip

\noindent{\bf Step Six: The proof of (\ref{511e1}) and (\ref{219e7}) }

\medskip

Let $y_1,y_2$ be two points on the same circle centered at $0$: $|y_1|=|y_2|$. Then the Green's representation formula gives
\begin{eqnarray}\label{223e9}
&& v_i^k(y_1)-v_i^k(y_2)\\
&=&\int_{B(0,\delta_k^{-1})}(G(y_1,\eta)-G(y_2,\eta))(\sum_{j=1}^na_{ij}h_j^k(\delta_k\eta)e^{v_j^k})d\eta \nonumber \\
&&+\hat h_i^k(y_1)-\hat h_i^k(y_2) \nonumber
\end{eqnarray}
where $\hat h^k=(\hat h_1^k,..,\hat h_n^k)$ is harmonic defined by $\hat h_i^k=v_i^k$ on $\partial B(0,\delta_k^{-1})$. Clearly
$$|\hat h_i^k(y)-\hat h_i^k(y')|=O(1),\quad \forall y,y'\in B(0,\delta_k^{-1})$$
because $v_i^k$ has bounded oscillation on $\partial B(0, \delta_k^{-1})$. Next
we claim
\begin{equation}\label{223e8}
|v_i^k(y_1)-v_i^k(y_2)|\le C,\quad \forall r\in (0,\delta_k^{-1}), \quad \forall |y_1|=|y_2|=r,\quad i\in I.
\end{equation}

We only need to verify (\ref{223e8}) for $r<\delta_k^{-1}/2$, as the case for $r>\frac 12\delta_k^{-1}$ is obvious. To evaluate the
expression for $v_i^k(y_1)-v_i^k(y_2)$, we first have
$$G(y_1,\eta)-G(y_2,\eta)=\frac 1{2\pi}\log \frac{|y_2-\eta |}{|y_1-\eta |}
+\frac 1{2\pi}\log \frac{|1-\delta_k^2\bar y_1\eta |}{|1-\delta_k^2\bar y_2\eta |} $$
where $\bar y_1$ is the conjugate of $y_1$ when it is considered as a complex number.
To prove (\ref{223e8}) we decompose $B(0,\delta_k^{-1})$ into the following four sets:
$$E_1:=\{\eta; \, |\eta |\leq r/2\}, \quad E_2:=\{\eta; \, |\eta-y_1|<|\eta -y_2|, \frac{r}{2}\leq |\eta |\leq 2r\}, $$
$$E_3:=\{\eta; \, |\eta -y_2|\le |\eta -y_1|, \, \frac r2\leq |\eta |\leq 2r \}, \quad E_4=B(0,\delta_k^{-1})\setminus (\cup_{i=1}^3E_i). $$
By (\ref{223e6})
$$e^{v_j^k(y)}\le (1+|y|)^{-2-\delta'},\quad j\in I, \quad |y|<\delta_k^{-1} $$
where $\delta'>0$. Applying this decaying estimate, we can prove the integral of the right hand side of (\ref{223e9}) over each $E_j$ is bounded.
We omit the proof because it is a standard computation. Thus (\ref{223e8}) is proved. Once (\ref{223e8}) is established, we have
\begin{equation}\label{223e10}
v_i^k(y)=\bar v_i^k(r)+O(1),\quad i\in I, \quad |y|=r.
\end{equation}
For $\bar v_i^k(r)$ we have
\begin{equation}\label{223e12}
(\bar v_i^k)'(r)=-\frac 1{2\pi r}\int_{B_r}\sum_{j=1}^na_{ij}h_j^k(\delta_ky)e^{v_j^k}dy.
\end{equation}
For $r>R$, the decay rates of $v_i^k$ in (\ref{223e6}) and (\ref{223e7}) imply
\begin{eqnarray}\label{223e11}
&&\frac 1{2\pi r}\int_{B_r}\sum_{j=1}^na_{ij}h_i^k(\delta_ky)e^{v_j^k}=\frac{1}{2\pi r}\bigg (\int_{B(0,\delta_k^{-1})}-\int_{B(0,\delta_k^{-1})\setminus
B_r} \bigg )\\
&=&\frac{m_i^k}r+O(r^{-1-\delta}).\nonumber
\end{eqnarray}
Using (\ref{223e11}) in (\ref{223e12}) we have
\begin{equation}\label{223e13}
\bar v_i^k(r)=v_i^k(0)-m_i^k\log (1+r)+O(1), \quad 0\le r<\delta_k^{-1}.
\end{equation}

Clearly (\ref{511e1}) and (\ref{219e7}) follow from (\ref{223e10}) and (\ref{223e13}).
This completes the proof of Proposition \ref{blowuplocal}.
 $\Box$

\begin{rem}
It is easy to see from \eqref{equofvi} that the sequence $u_i^k(0)-M_k$ is bounded for $1\le i\le l$ and
tends to $-\infty$ for $l+1\le i\le n$.
\end{rem}

\section{The strong interaction between bubbles}

In this section, we suppose $u^k$ has two blow-up points $p_1$ and $p_2$.
By Proposition \ref{blowuplocal}, at each blow-up point $p_i$, $u^k$ after scaling will converges to an entire solution of a subsystem \eqref{equofvi}.
The following question naturally arises:
\begin{center}
Are these two entire solutions the same?
\end{center}
The following theorem will answer this question affirmatively.

\begin{prop}\label{bubinter} Let $\Omega_0$ be an open and bounded set with smooth boundary, $p_1,p_2 \in \Omega_0$ be two distinct points.
Suppose $u^k=(u_1^k,..,u_n^k)$ satisfies (\ref{uikeq}), (\ref{219e2}) and (\ref{219e4}) on $\Omega_0$ and and $h^k$ satisfies (\ref{hika}) over
$\Omega_0$ as well. Let $p_1,p_2$ be the only two blow-up points of $u^k$ on $\Omega$:
$$\exists  x_{tk}\to p_t, i_t\in I, \quad \mbox{ such that }\,\, u_{i_t}^k(x_{tk})\to \infty, \quad t=1,2. $$
$$\max_K u_i^k\le C(K),\quad \forall K\subset \subset \Omega_0\setminus \{p_1,p_2\},\quad i\in I. $$
Then for $\delta<\frac 12|p_1-p_2|$
$$\lim_{k\to \infty}\int_{B(p_1,\delta)}h_i^ke^{u_i^k}dx=
\lim_{k\to \infty}\int_{B(p_2,\delta)}h_i^ke^{u_i^k}dx,\quad i\in I. $$
\end{prop}

\noindent{\bf Proof of Proposition \ref{bubinter}:}  Let
$$\sigma_i^k=\frac{1}{2\pi}\int_{B(p_1,\delta)}h_i^ke^{u_i^k}dx,
\quad \bar \sigma_i^k=\frac{1}{2\pi}\int_{B(p_2,\delta)}h_i^ke^{u_i^k}dx,$$
for $\delta>0$ small and $i\in I$. Also we let
$$m_i^k=\sum_{j=1}^na_{ij}\sigma_j^k,\quad \bar m_i^k=\sum_{j=1}^na_{ij}\bar \sigma_j^k,\quad i\in I. $$
We use $\sigma_i, m_i, \bar \sigma_i$ and $\bar m_i$ to denote the limit of $\sigma_i^k,m_i^k,\bar \sigma_i^k$ and $\bar m_i^k$, respectively.
Let
$$M_k=\max_{i\in I}u_i^k(x),\quad x\in B(p_1,\delta),$$
and $M_k$ is attained by some component of $u^k$ at $p_{1k}$ which tends to $p_1$. $\bar M_k$ and $p_{2k}$ can be defined similarly.
By comparing the value of $u_i^k$ over $\Omega_0\setminus (B(p_1,\delta)\cup B(p_2,\delta))$, using Proposition \ref{blowuplocal} we have
\begin{equation}\label{compeq}
\frac{m_i^k-2}2M_k+(M_k-u_i^k(p_{1k}))=\frac{\bar m_i^k-2}2\bar M_k+(\bar M_k-u_i^k(p_{2k}))+O(1).
\end{equation}
Here we remind the reader that, for example around $p_1$, if the first $l$ components of $u^k$ converge to a system of $l$ equations after scaling, then
$M_k-u_i^k(p_{1k})$ are uniformly bounded for $1\le i\le l$. In this case, $M_k-u_i^k(p_{1k})$ can be combined with the $O(1)$ term.
For $i>l$, $M_k-u_i^k(p_{ik})$ tends to $+\infty$. The right hand side of (\ref{compeq}) can also be understood this way.
For each $i\in I$, if
$$M_k-u_i^k(p_{1k})>\bar M_k-u_i^k(p_{2k})$$
we let
$$l_i^k=(M_k-u_i^k(p_{1k}))-(\bar M_k-u_i^k(p_{2k})),\quad \bar l_i^k=0. $$
On the other hand if
$$M_k-u_i^k(p_{1k})\le \bar M_k-u_i^k(p_{2k})$$
we let
$$l_i^k=0,\quad \bar l_i^k=(\bar M_k-u_i^k(p_{2k}))-(M_k-u_i^k(p_{1k})). $$
Set
$$I_1:=\{i\in I;\,\, \lim_{k\to \infty}\frac{l_i^k}{\bar M_k}>0.\, \},\quad
I_2:=\{i\in I;\,\, \lim_{k\to \infty}\frac{\bar l_i^k}{\bar M_k}>0.\,\}$$
and $I_3$ be the compliment of $I_1\cup I_2$. From this definition we see easily that $I_1\cap I_2=\emptyset$.

We claim that $I_1$ is empty. We prove this by contradiction. Suppose $I_1$ is not empty, then we consider the following two cases:
$I_2$ is not empty or $I_2$ is empty.

\medskip

\noindent{\bf Case one: $I_2\neq \emptyset$}

\medskip

Let
$$\lambda=\lim_{k\to \infty}\frac{M_k}{\bar M_k},\quad \delta_i=\lim_{k\to \infty}\frac{l_i^k}{\bar M_k},\quad
\bar \delta_i=\lim_{k\to \infty}\frac{\bar l_i^k}{\bar M_k}. $$
We claim that all these limits exist. Indeed, using the definitions of $l_i^k$ and $\bar l_i^k$, (\ref{compeq}) can be
written as
$$\frac{m_i^k-2}2\cdot \frac{M_k}{\bar M_k}+\frac{l_i^k}{\bar M_k}=
\frac{\bar m_i^k-2}2+\frac{\bar l_i^k}{\bar M_k}+\circ(1). $$
Take $i\in I_1$, the RHS tends to $(\bar m_i-2)/2$, which implies that along a subsequence, the two terms on the LHS are
$\frac{m_i-2}2\lambda$ and $\delta_i$. On the other hand, take $j\in I_2$, the LHS tends to $\frac{m_j-2}2\lambda$, then the RHS has to tend to $\frac{\bar m_j-2}2+\bar \delta_j$.
Now (\ref{compeq})
can be written as
\begin{equation}\label{compeq1}
\lambda \frac{m_i-2}2+\delta_i=\frac{\bar m_i-2}2+\bar \delta_i,\quad i\in I.
\end{equation}
From the definition of $\delta_i$ and $\bar \delta_i$, we observe that for each $i\in I_1$, $\bar \delta_i=0$ and for $i\in I_2$, $\delta_i=0$.
By (\ref{511e1}) and (\ref{219e7}) of Proposition \ref{blowuplocal}, we have
\begin{equation}\label{subtleob}
\sigma_i=0,\,\, i\in I_1; \qquad \bar \sigma_i=0,\,\,i\in I_2.
\end{equation}
Since $\delta_i=0$ for $i\not\in I_1$ and $\sigma_i=0$ for $i\in I_1$, we have $\sigma_i\delta_i=0$ for all $i\in I$.
Similarly, $\overline \sigma_i \bar \delta_i=0$ for all $i\in I$.

Without loss of generality we assume
$$I_1=\{1\le i\le i_0\}.$$
For each $i\in I_2$, the fact $\bar \sigma_i=0$ yields
\begin{equation}\label{513e2}
0=\bar \sigma_i=\sum_{j\in I}a^{ij}\bar m_j=\sum_{j\in I_2}a^{ij}\bar m_j+\sum_{j\not \in I_2}a^{ij}\bar m_j.
\end{equation}
We observe that the last term is positive, because $\bar m_j>2$ and there exists $a^{ij}> 0$ for some $i\in I_2$ and $j\not \in I_2$. Multiplying
 $\bar \delta_i$ to the last term and taking the summation on $i$ for all $i$ in $I_2$, we have
\begin{equation}\label{513e1}
\sum_{i\in I_2}(\sum_{j\not \in I_2}a^{ij}\bar m_j)\bar \delta_i>0.
\end{equation}
Combining (\ref{513e2}) and (\ref{513e1}), we have
$$\sum_{i,j\in I_2}a^{ij}\bar m_i\bar \delta_j<0. $$
Trivially, there exists $\tilde i\in I_2$ such that
\begin{equation}\label{513e3}
\sum_{j\in I_2}a^{\tilde ij}\bar \delta_j<0.
\end{equation}
Multiplying $a^{\tilde ij}$ to both sides of (\ref{compeq1}) (with $i$ replaced by $j$) and taking the summation on $j$, it leads to
\begin{equation}\label{513e4}
\sum_{j\in I}a^{\tilde ij}(\frac{m_j-2}2)\lambda+\sum_{j\in I}a^{\tilde ij}\delta_j
=\sum_{j\in I}a^{\tilde ij}(\frac{\bar m_j-2}2)+\sum_{j\in I}a^{\tilde ij}\bar \delta_j.
\end{equation}
Using the definition of $\sigma_{\tilde i}$ as well as $\delta_i=0$ for $i\not \in I_1$, we can write the LHS of (\ref{513e4}) as
$$\frac 12\lambda\sigma_{\tilde i}-\lambda\sum_{j\in I}a^{\tilde ij}+\sum_{j\in I_1}a^{\tilde ij}\delta_j. $$
Since the first term and the third term are both nonnegative, the LHS is no less than the second term. Similarly the RHS can be written as
$$\frac 12\bar \sigma_{\tilde i}-\sum_{j\in I}a^{\tilde ij}+\sum_{j\in I_2}a^{\tilde ij}\bar \delta_j. $$
Note that we have used $\bar \delta_i=0$ for $i\not \in I_2$. The first term of the above is $0$ (because $\tilde i\in I_2$) and the last term is negative (because of (\ref{513e3})). Therefore the RHS is strictly less than the second term. Putting the estimates on both sides together we have
$$-\lambda \sum_{j\in I}a^{\tilde ij}<-\sum_{j\in I}a^{\tilde ij}. $$
Since $\sum_{j\in I}a^{\tilde ij}\ge 0$ (({\bf H2})) we conclude $\lambda>1$. On the other hand, by exchanging $I_1$ and $I_2$ in the above argument, we obtain
$\lambda<1$. Thus we have ruled out the first case.

\medskip

\noindent{\bf Case two: $I_2=\emptyset$}

\medskip

One immediately has $\bar \delta_i=0$ for all $i\in I$. Hence, the limits $\lambda=\lim_{k\to \infty}M_k/\bar M_k$ and
$\delta_k=\lim_{k\to \infty}l_i^k/\bar M_k$ both exist and (\ref{compeq1}) holds with $\bar \delta_k=0$.
Here we recall that both $\sigma=(\sigma_1,..,\sigma_n)$ and $\bar \sigma=(\bar \sigma_1,..,\bar \sigma_n)$ satisfy
$$4\sum_{i\in I}\sigma_i=\sum_{i,j\in I}a_{ij}\sigma_i\sigma_j$$
which can be written as
\begin{equation}\label{513e5}
\sum_{i,j\in I}a^{ij}(\frac{m_i-2}2)(\frac{m_j-2}2)=\sum_{i,j\in I}a^{ij}.
\end{equation}
We further remark that
\begin{equation}\label{513e6a}
\sum_{i,j\in I}a^{ij}>0
\end{equation}
because for each $i$, $\sum_{j\in I}a^{ij}\ge 0$ and $A^{-1}$ is non-singular. To prove our result, we need another fact:
\begin{equation}\label{513e6}
\mbox{All the eigenvalues of ${\bf F}$ are nonpositive},
\end{equation}
where
$${\bf F}=(a^{ij})_{i_0\times i_0}\quad i,j\in I_1=\{1,..,i_0\}. $$
Indeed, let $\mu$ be the largest eigenvalue of ${\bf F}$ and $\eta=(\eta_1,..,\eta_{i_0})$ be an eigenvector corresponding to $\mu$. Here
 $\eta$ is the vector that attains
 $$\max_{{\bf v}\in \mathbb R^n} {\bf v^T\, F \, v},\quad {\bf v^Tv}=1. $$
Since $a^{ij}\ge 0$ for all $i\neq j$, we can choose $\eta_i\ge 0$ for all $i\in I_1$.
For each $i\in I_1$,
$$ 0=\sigma_i=\sum_{j\in I_1}a^{ij}m_j+\sum_{j\not \in I_1}a^{ij}m_j. $$
Plainly by ({\bf H2})
$$\sum_{j\in I_1}a^{ij}m_j\le 0,\quad i\in I_1. $$
Multiplying $\eta_i$ on both sides and taking the summation on $i$, then we have
$$0\ge \sum_{i,j\in I_1}a^{ij}\eta_im_j=\sum_{j\in I_1}\mu\eta_jm_j. $$
Since each $\eta_i\ge 0$ (one of them is strictly positive) and $m_i>0$ for $i\in I_1$, we have $\mu\le  0$, and it proves \eqref{513e6}.
Now we go back to our proof to rule out case two.

Since $\bar \delta_i=0$ in (\ref{compeq1}), the Pohozaev identity \eqref{513e5} for $\bar \sigma$ can be written as
$$\sum_{i,j\in I}a^{ij}(\frac{m_i-2}2\lambda+\delta_i)(\frac{m_j-2}2\lambda+\delta_j)=\sum_{ij}a^{ij}. $$
Expanding the LHS of the above and using (\ref{513e5}) again for $\sigma_i$, we obtain
\begin{equation}\label{513e7}
\lambda^2\sum_{i,j\in I}a^{ij}+2\lambda\sum_{i,j\in I}a^{ij}(\frac{m_i-2}2)\delta_j+\sum_{i,j\in I_1}a^{ij}\delta_i\delta_j
=\sum_{i,j}a^{ij}.
\end{equation}
The third term of LHS is nonpositive by (\ref{513e6}). The second term of the LHS can be written as
$$\lambda \sum_{j\in I}(\sigma_j-2\sum_{i\in I}a^{ij})\delta_j=-\lambda\sum_{j\in I}(\sum_{i\in I}a^{ij})\delta_j\leq 0 $$
because $\sigma_j \delta_j=0$ for all $j\in I$.
Thus we conclude from (\ref{513e7}) that $\lambda\ge 1$.
On the other hand from $\sigma_i=0$ ($i\in I_1$), argued as \eqref{513e3} we obtain an index $\tilde i\in I_1$ such that
\begin{equation}\label{513e8}
\sum_{j\in I_1}a^{\tilde ij}\delta_j<0.
\end{equation}
Then as we did for \eqref{513e4}, we obtain
$$\lambda\sum_{j\in I}a^{\tilde ij}(\frac{m_j-2}2)+\sum_{j\in I}a^{\tilde ij}\delta_j=\sum_{j\in I}a^{\tilde ij}(\frac{\bar m_j-2}2). $$
Following the same calculation as before we obtain
$$0>\sum_{j\in I_1}a^{\tilde ij}\delta_j=\frac{\bar \sigma_{\tilde i}}2+(\lambda-1)\sum_{j\in I}a^{\tilde ij}, $$
which forces $\lambda$ to be less than $1$. Thus, we have obtained a contradiction to $\lambda\ge 1$. Case two is also ruled out. Thus we have proved that $I_1$ has to be empty. Using exactly the same argument we
also have $I_2=\emptyset$.

\medskip

Since $I_1=I_2=\emptyset$, (\ref{compeq1}) becomes
$$\lambda \frac{m_i-2}2=\frac{\bar m_i-2}2,\quad i\in I. $$
Using (\ref{513e5}) for both $(m_1,..,m_n)$ and $(\bar m_1,..,\bar m_n)$ we have $\lambda=1$. Consequently $\sigma_i=\bar \sigma_i$ for all $i\in I$. Proposition \ref{bubinter} is established. $\Box$

\section{Proof of Theorem \ref{main1} and Theorem \ref{main2}}

Let $u=(u_1,..,u_n)$ be a solution of (\ref{mainsys}). Define
\begin{equation}\label{64e1}
v_i=u_i-\log \int_Mh_ie^{u_i}dV_g.
\end{equation}
Clearly $v=(v_1,..,v_n)$ satisfies
\begin{equation}\label{64e3}
\int_Mh_ie^{v_i}dV_g=1
\end{equation}
and
\begin{equation}\label{eq51a}
\Delta_g v_i+\sum_{j=1}^n\rho_j a_{ij}(h_je^{v_j}-1)=0,\quad i\in I.
\end{equation}
To prove the a priori bound for $u$, we only need to establish
\begin{equation}\label{vbound}
|v_i(x)|\le C,\quad \forall x\in M,\quad i\in I.
\end{equation}
Indeed, once we have (\ref{vbound}), for $u$ we have
\begin{equation}\label{ubound}
\log \int_Mh_ie^{u_i}-C\le u_i(x)\le \log \int_Mh_ie^{u_i}+C,\quad \forall x\in M.
\end{equation}
Since $u\in$  $\hdot^1(M)$, there exists $x_0$ such that $u(x_0)=0$. Using this in (\ref{ubound}) we have
\begin{equation}\label{64e2}
-C\le \log \int_Mh_ie^{u_i}\le C.
\end{equation}
With (\ref{64e1}) and (\ref{64e2}) we see that a bound for $u$ can be obtained from the bound for $v$. To prove (\ref{vbound}) we only need to prove
\begin{equation}\label{vbound1}
v_i\le C,\quad i\in I.
\end{equation}
because a lower bound for $v_i$ can be obtained easily from the upper bound in (\ref{vbound1}) by standard elliptic estimate. So we only need to establish (\ref{vbound1}).

We prove (\ref{vbound1}) by contradiction. Suppose there are solutions $v^k$ to (\ref{eq51a})
such that $\max_{M, i\in I}v_i^k(x)\rightarrow+\infty$. We consider the following two cases.

 \medskip

\noindent{\bf Case one: } $\rho_i^k\rightarrow \rho_i>0$ as $k\rightarrow+\infty$ for all $i\in I$.

The equation for $v^k$ is
\begin{equation}\label{eq51}
\Delta_g v_i^k+\sum_{j=1}^n\rho_j^k a_{ij}(h_je^{v_j^k}-1)=0,\quad i\in I.
\end{equation}
In \cite{linzhang1} the authors prove a Brezis-Merle type lemma (Lemma 4.1) which guarantees that there exists a positive constant $\epsilon_0>0$ such that if
$$\int_{B(p,r_0)}e^{v_i^k}dx\leq \epsilon_0 \quad \mbox{ for all } \quad i\in I, $$ then
\begin{equation}
v_i^k(x)\leq C,\quad x\in \overline B\left(p,\frac{r_0}{2}\right).
\end{equation}

Thus, $v^k$ blows up only at a finite set $\{p_1,\cdots, p_N\}$. Since $v_i^k(x)$ is uniformly bounded from above
in any compact set of $M\setminus\{p_1,\cdots,p_N\}$, by \eqref{eq51}, $v_i^k$ converges to $\sum_{l=1}^N m_i(p_l)G(x,p_l)$ in
$C_{loc}^\infty(M\setminus\{p_1,\cdots,p_N\})$, where
\begin{equation}
\begin{cases}
&m_i(p_l)=\sum_{j\in I}a_{ij}\sigma_j(p_l),\\
&\sigma_j(p_l)=\lim_{k\to \infty}\frac{1}{2\pi}\int_{B(p_l,\delta_0)}\rho_j^kh_j e^{v_j^k}dV_g
\end{cases}
\end{equation}
for some $\delta_0>0$ such that $B(p_l, 2\delta_0)\cap B(p_{l'}, 2\delta_0)=\emptyset$, $l\ne l'$. Here, $G(x,p)$ is the Green function:
\begin{equation*}
\begin{cases}
&-\Delta_{g,x}G(x,p)=\delta_p-1,\\
&\int_M G(x,p) dV_g(x)=0.
\end{cases}
\end{equation*}

To apply Proposition \ref{bubinter}, we rewrite \eqref{eq51} in local coordinates.
For $p\in M$, let $y=(y^1, y^2)$ be the isothermal coordinates near $p$ such that $y_p(p)=(0,0)$ and $y_p$ depends smoothly on $p$.
In this coordinate, $ds^2$ has the form
\begin{equation*}
e^{\phi(y_p)}\left[(dy^1)^2+(dy^2)^2\right],
\end{equation*}
where $\nabla \phi(0)=0$, $\phi(0)=0$. Also near $p$ we have
\begin{equation*}
\Delta_{y_p}\phi=-2 K e^{\phi},\quad\text{where \,}K\text{ is the Gauss curvature}.
\end{equation*}
When there is no ambiguity, we write $y=y_p$ for simplicity. In this local coordinate, \eqref{eq51} is of the form:
\begin{equation}\label{uikequation}
-\Delta v_i^k=e^{\phi}\sum_{j=1}^n a_{ij}\rho_j^k \left(
h_j e^{v_j^k}-1
\right)\quad\text{in \,}B(0, \delta_0),\quad i\in I.
\end{equation}
Let $f_i^k$ be defined as
\begin{equation*}
-\Delta f_i^k=-e^{\phi}\sum_{j=1}^n \rho_j^k a_{ij}\quad\text{in \,}B(0,\delta_0),\quad i\in I,
\end{equation*}
and $f_i^k(0)=|\nabla f_i^k(0)|=0$. Let $\tilde v_i^k=v_i^k-f_i^k$ and
\begin{equation*}
H_i^k= e^{\phi}\rho_i^k e^{f_i^k} h_i.
\end{equation*}
Then \eqref{uikequation} becomes
\begin{equation}\label{616e1}
-\Delta \tilde v_i^k=\sum_{j=1}^n a_{ij} H_j^k e^{\tilde v_j^k}\quad\text{in \,}B(0,\delta_0).
\end{equation}
Here, we observe that $\int_{B(0,\delta_0)}H_i^k e^{\tilde u_i^k}dx=\int_{B(0,\delta_0)}\rho_i^kh_i e^{v_i^k}dV_g$.

Since $v_i^k$ converges in $M\setminus\bigcup_{j=1}^N B(p_j, 2\delta_0)$, we have
\begin{equation}
|\tilde v_i^k(x)-\tilde v_i^k(y)|\leq C,\quad \forall x,y\in M\setminus\bigcup_{j=1}^N B(p_j,2\delta_0),\quad i\in I.
\end{equation}
By \eqref{511e1} of Proposition \ref{blowuplocal}, we also have
\begin{equation}\label{cond1}
\int_{M\setminus \bigcup_{j=1}^N B(p_j,\delta_0)}h_i e^{v_i^k}dV_g\rightarrow 0,\quad i\in I,
\end{equation}
and by Proposition \ref{bubinter},
\begin{equation}\label{cond2}
\lim_{k\rightarrow\infty}\int_{B(p_l,\delta_0)}\rho_i^k h_i e^{v_i^k}dV_g
=\lim_{k\rightarrow\infty}\int_{B(p_m,\delta_0)}\rho_i^kh_i e^{v_i^k}dV_g
\end{equation}
for $i\in I$ and for any pair of integers $l, m$ between $1$ and $N$. \eqref{cond1} combined with \eqref{cond2} yields for $i\in I$ and $j\in\{1,2,\cdots,N\}$,
\begin{equation*}
\sigma_i=\lim_{k\rightarrow\infty}\frac{1}{2\pi}\int_{B(0,\delta_0)}H_i^k e^{\tilde v_i^k}dx=\lim_{k\rightarrow\infty}\frac{1}{2\pi}\int_{B(p_j,\delta_0)} \rho_i^k h_i e^{v_i^k} dV_g=\frac{\rho_i}{2\pi N}.
\end{equation*}
On the other hand, $(\sigma_1,\cdots,\sigma_n)$ satisfies the Pohozaev identity:
\begin{equation}
4\sum_{i\in I}\sigma_i=\sum_{i,j\in I}a_{ij} \sigma_i\sigma_j.
\end{equation}
Consequently,
\begin{equation*}
8\pi N \sum_{i=1}^n \rho_i=\sum_{i,j=1}^n a_{ij}\rho_i\rho_j.
\end{equation*}
Thus, a contradiction to the assumption of the theorem.

\medskip

\noindent{\bf Case two:} $\lim_{k\to \infty}\rho_i^k=\rho_i>0$, $i=1,..,l$, $\lim_{\to \infty}\rho_i^k=0$ for $i>l$.

\medskip

Let $M_k=\max\{ v_1^k,..,v_l^k\}$ and $\bar M_k=\{v_{l+1}^k,..,v_n^k\}$. We first show
\begin{equation}\label{63e1}
\bar M_k-M_k\le C
\end{equation}
by contradiction. Suppose $\bar M_k-M_k\to \infty$, let
$$V_i^k(y)=v_i^k(e^{-\frac{\bar M_k}2}y+p_k)-\bar M_k$$
where $p_k$ is where $\bar M_k$ is attained: $v_{i_0}^k(p_k)=\bar M_k$. Clearly $i_0>l$. Thanks to the fact that $V_i^k\to -\infty$ for $i\le l$ and $\rho_i^k\to 0$ for $i>l$, $V_{i_0}^k$ converges in $C^2_{loc}(\mathbb R^2)$ to
$$\left\{\begin{array}{ll}
-\Delta V_{i_0}=0,\quad \mathbb R^2,\\
V_{i_0}(0)=0,\quad V_{i_0}\le 0.
\end{array}
\right.
$$
Clearly $V_{i_0}\equiv 0$, $\int_{B_R}e^{V_{i_0}}$ can be arbitrarily large if $R$ is large, this is a contradiction to (\ref{64e3}). (\ref{63e1}) is proved.

We use the same notations as in case one. Let $p_1,..,p_N$ be blowup points for $v_i^k$. Then around each blowup point, say $p_1$, the equation for $v^k$ can be written in local coordinates as (\ref{616e1}) with $\tilde v_i^k$ and $H_i^k$ defined the same as in case one. Without loss of generality we assume $\rho_i^k>0$ for all $k$ and $l+1\le i\le L$ and $\rho_i^k=0$ for all $k$ and all $L+1\le i\le n$. Then we observe from the definition of $H_i^k$ that $H_i^k\to 0$ for $l+1\le i\le L$ and $H_i^k=0$ for $i>L$. To reduce case two to case one, we need to adjust the terms involved with these vanishing $H_i^k$s. To do this we set $\hat f_i^k$ as
$$\left\{\begin{array}{ll}
-\Delta \hat f_i^k=\sum_{j=L+1}^na_{ij}e^{\tilde v_j^k-M_k},\quad B(0,\delta),\\
\\
\hat f_i^k=0\quad \mbox{ on }\quad \partial B(0,\delta).
\end{array}
\right.
$$
Since $\max v_i^k-M_k$ is uniformly bounded for all $i$, we have
$$\|\hat f_i^k\|_{C^1}\le C $$
for some $C$ independent of $k$. Now we define
$$\hat v_i^k=\left\{\begin{array}{ll}
\tilde v_i^k+\hat f_i^k,\quad i=1,..,l,\\
\tilde v_i^k+\log \rho_i^k+\hat f_i^k,\quad l+1\le i\le L,\\
\tilde v_i^k-M_k+\hat f_i^k,\quad L+1\le i\le n.
\end{array}
\right.
$$
and
$$\hat H_i^k=\left\{\begin{array}{ll}
H_i^ke^{-\hat f_i^k},\quad 1\le i\le l,\\
\\
\frac{H_i^k}{\rho_i^k}e^{-\hat f_i^k}=e^{\phi+f_i^k-\hat f_i^k}h_i,\quad l+1\le i\le L, \\
\\
e^{-\hat f_i^k},\quad L+1\le i\le n.
\end{array}
\right.
$$
Easy to see there exists $c>0$ independent of $k$ such that
$$\frac 1c\le \hat H_i^k\le c, \quad |\nabla \hat H_i^k|\le c, \quad B(0,\delta). $$
On the other hand $\hat v_i^k$ satisfies
$$-\Delta \hat v_i^k=\sum_{j\in I}a_{ij}\hat H_j^k e^{\hat v_j^k},\quad B(0, \delta), \quad i\in I. $$
Easy to observe that $\max \hat v_i^k-M_k\to -\infty$ for $i=l+1,...,n$. Therefore case two is reduced to case one, which gives
$$\sigma_i(p_t)=\sigma_i(p_m)\quad \forall t, m \in \{1,..,N\}, \quad \forall 1\le i\le l. $$
Note that $\sigma_i(p_t)=0$ for all $i>l$ and all $t$ because $\rho_i^k\to 0$ for $i>l$. Then as in case one we obtain a contradiction.
Theorem \ref{main1} is established. $\Box$

\bigskip

\noindent{\bf Proof of Theorem \ref{main2}:}

\medskip

Theorem \ref{main2} will be discussed in two cases.

\medskip

{\bf Case 1.} One of $a_{ii}$ is positive.

\medskip

We may suppose $a_{11}>0$.
Thanks to Theorem \ref{main1} the Leray-Schauder degree of \eqref{mainsys} for $\rho\in \mathcal{O}_{N-1}$ is equal to the degree for
the following specific system corresponding to $(\rho_1,0,..,0)$:

\begin{equation}\label{mainsystrans}
\Delta_g u_1+\rho_1 a_{11}\left(\frac{h_1 e^{u_1}}{\int_M h_1 e^{u_1} dV_g}-1\right)=0,
\end{equation}
$$
\Delta_g u_j+\rho_1 a_{j1}\left(\frac{h_1 e^{u_1}}{\int_M h_1 e^{u_1}dV_g}-1\right)=0\quad\text{for \,}j\geq 2.
$$
where $\rho_1$ satisfies
\begin{equation}\label{thm13proof1}
8\pi (N-1)<a_{11}\rho_1<8\pi N.
\end{equation}
Easy to see $(\rho_1,0,..,0)\in \mathcal{O}_{N-1}$.
Using Theorem 1.2 of \cite{chenlin2}, we obtain the degree counting formulas \eqref{dcount} in this case.

\medskip

\noindent
{\bf Case 2.} $a_{ii}=0$ for all $i\in I$.

\medskip

By Lemma \ref{irredu}, $a_{12}>0$. The degree counting formula of (\ref{mainsys}) for $\rho\in \mathcal{O}_N$ can be computed by
the degree of the following specific system
\begin{equation}\label{equu1u2}
\begin{cases}
&\Delta_g u_1+a_{12}\rho_2\left(\frac{h_2 e^{u_2}}{\int_M h_2e^{u_2}dV_g}-1\right)=0,\\
&\Delta_g u_2+a_{12}\rho_1\left(\frac{h_1 e^{u_1}}{\int_M h_1 e^{u_1}dV_g}-1\right)=0,\\
&\Delta_g u_i+\rho_1a_{i1}\left(\frac{h_1 e^{u_1}}{\int_M h_1 e^{u_1}dV_g}-1\right)\\
&\qquad +
\rho_2a_{i2}\left(\frac{h_2 e^{u_2}}{\int_M h_2e^{u_2}dV_g}-1\right)=0,\quad i\ge 3.
\end{cases}
\end{equation}
where $\rho_1,\rho_2$ satisfy
\begin{equation}\label{assump513}
8\pi (N-1)(\rho_1+\rho_2)<2 a_{12}\rho_1\rho_2<8\pi N(\rho_1+\rho_2).
\end{equation}
Easy to see $(\rho_1,\rho_2,0,..,0)\in \mathcal{O}_{N-1}$.
Now we consider the special case of \eqref{equu1u2}: $\rho_1=\rho_2$ and $h_1=h_2=h$. In this case, the maximum principle implies $u_1=u_2+c$
for some constant $c$. Since $u_1,u_2 $ are both in $\hdot^1(M)$, $c=0$. Then the first two equations of \eqref{equu1u2} turn out to be:
\begin{equation}\label{lssingleeq}
\Delta_g u+a_{12}\rho\left(\frac{h e^u}{\int_M h e^u dV_g}-1\right)=0,
\end{equation}
where $\rho$ satisfies
\begin{equation*}
8\pi (N-1)<a_{12}\rho<8\pi N.
\end{equation*}
Hence, again the Leray-Schauder degree for equation \eqref{mainsys} can be reduced to the Leray-Schauder degree for the single equation \eqref{lssingleeq}.
By applying Theorem 1.2 in \cite{chenlin2}, the degree counting formulas \eqref{dcount} is also obtained in this case. This completes the proof of Theorem \ref{main2}. $\Box$

\section{Proof of Theorem \ref{dirthm}}

For equation \eqref{systemtrans}, we have to show $u^k$ never blows up near the boundary $\partial\Omega$. This fact is standard, we include the argument for the convenience of the reader (see \cite{mawei}). Since $\Omega$ is a bounded set with smooth boundary, there is a uniform constant $r_0$ such that for any point on $\partial \Omega$, there is a ball of radius $r_0$ tangent to $\partial \Omega$ at this point from the outside. Let $x_0\in \partial \Omega$ and $B(x_1,\lambda)$ be a ball tangent to $\partial \Omega$ at $x_0$ from the outside. $\lambda\le r_0$ will be determined later. Let
$$H_i=\frac{\rho_ih_i}{\int_{\Omega}h_ie^{u_i}dx},\quad \rho_i>0, \,\, i\in I$$
then the equation for $u_i$ becomes
$$-\Delta u_i=\sum_{j=1}^na_{ij}H_je^{u_j},\quad \Omega,\quad i\in I.  $$
For $H_i$ we obviously have
\begin{equation}\label{Higrad}
|\nabla \log H_i(x)|\le C,\quad \forall x\in \Omega.
\end{equation}
Let $y=x-x_1$ and
$$u_i^{\lambda}(y)=u_i(x_1+\lambda^2\frac{y}{|y|^2}),\quad i\in I. $$
Then
$$-\Delta u_i^{\lambda}(y)=\sum_{j=1}^na_{ij}\bigg (\frac{\lambda^4}{|y|^4} H_j(x_1+\lambda^2\frac{y}{|y|^2}) \bigg )e^{u_j^{\lambda}(y)} \quad \mbox{in}\quad
\Omega^{\lambda} $$
where $\Omega^{\lambda}$ is the image of $\Omega$ under the Kelvin transformation. Moreover, we have
$$u_i^{\lambda}=0\quad \mbox{on}\quad \partial \Omega^{\lambda}. $$
Let
$$\bar H_i(y)=\frac{\lambda^4}{|y|^4} H_i(x_1+\lambda^2\frac{y}{|y|^2}).$$
Using (\ref{Higrad}), we see by direct computation that in a small neighborhood of $x_0$,
$\bar H_i$ is strictly decreasing in the outer normal direction to $\partial \Omega^{\lambda}$, as long as $\lambda$ is small. The smallness of
the neighborhood of $x_0$ and $\lambda$ can both be represented by $\epsilon_0$ which depends on the usual constants. Thus we have the monotonicity
of $\bar H_i$ in a neighborhood of the whole $\partial \Omega^{\lambda}$. Using the standard moving plane argument we see $u_i^{\lambda}$ is increasing along the inner normal of $\partial \Omega^{\lambda}$ in a small neighborhood of $\partial \Omega^{\lambda}$, which implies that for any sequence of function $u^k$ of (\ref{systemtrans}), no blowup point for $u^k$ exists in a fixed neighborhood of $\partial \Omega$.
Then the remaining part of the proof of Theorem \ref{dirthm} is the same as Theorem \ref{main1} and Theorem \ref{main2}.
So, the details are omitted here. $\Box$

\end{document}